\documentclass[english,12pt]{article}
\usepackage[T1]{fontenc}
\usepackage[latin9]{inputenc}
\usepackage{geometry}
\usepackage{babel}
\usepackage{mathrsfs}
\usepackage{bm}
\usepackage{amsmath}
\usepackage{amsthm}
\usepackage{amssymb}
\usepackage{stackrel}
\usepackage{rotfloat}
\usepackage{setspace}
\usepackage[authoryear]{natbib}
\PassOptionsToPackage{normalem}{ulem}
\usepackage{ulem}
\usepackage[unicode=true]
 {hyperref}

\makeatletter

\providecommand{\tabularnewline}{\\}

\newcommand{\lyxaddress}[1]{
	\par {\raggedright #1
	\vspace{1.4em}
	\noindent\par}
}
\theoremstyle{plain}
\newtheorem{lem}{\protect\lemmaname}
\theoremstyle{remark}
\newtheorem{rem}{\protect\remarkname}
\theoremstyle{plain}
\newtheorem{thm}{\protect\theoremname}
\newenvironment{lyxlist}[1]
	{\begin{list}{}
		{\settowidth{\labelwidth}{#1}
		 \setlength{\leftmargin}{\labelwidth}
		 \addtolength{\leftmargin}{\labelsep}
		 }}
	{\end{list}}

\makeatother

\providecommand{\lemmaname}{Lemma}
\providecommand{\remarkname}{Remark}
\providecommand{\theoremname}{Theorem}

\begin{document}
\title{PanIC: consistent information criteria for general model selection
problems}
\author{Hien Duy Nguyen}
\maketitle

\lyxaddress{School of Computing, Engineering and Mathematical Sciences, La Trobe
University, Bundoora, Australia. Institute of Mathematics for Industry,
Kyushu University, Fukuoka, Japan.\\
Email: \href{mailto:h.nguyen5@latrobe.edu.au}{h.nguyen5@latrobe.edu.au}.}
\begin{abstract}
Model selection is a ubiquitous problem that arises in the application
of many statistical and machine learning methods. In the likelihood
and related settings, it is typical to use the method of information
criteria (IC) to choose the most parsimonious among competing models
by penalizing the likelihood-based objective function. Theorems guaranteeing
the consistency of IC can often be difficult to verify and are often
specific and bespoke. We present a set of results that guarantee consistency
for a class of IC, which we call PanIC (from the Greek root '\emph{pan}',
meaning '\emph{of everything}'), with easily verifiable regularity
conditions. The PanIC are applicable in any loss-based learning problem
and are not exclusive to likelihood problems. We illustrate the verification
of regularity conditions for model selection problems regarding finite
mixture models, least absolute deviation and support vector regression,
and principal component analysis, and we demonstrate the effectiveness
of the PanIC for such problems via numerical simulations. Furthermore,
we present new sufficient conditions for the consistency of BIC-like
estimators and provide comparisons of the BIC to PanIC.
\end{abstract}

\section{\label{sec:Introduction}Introduction}

Let $\left(\Omega,\mathcal{F},\text{Pr}\right)$ be a probability
space with typical element $\omega$ and expectation operator $\text{E}$,
and suppose that $\bm{X}:\Omega\rightarrow\mathbb{X}\subset\mathbb{R}^{d}$
is a random variable on the image probability space $\left(\mathbb{X},\mathcal{B}\left(\mathbb{X}\right),\text{Pr}_{\bm{X}}\right)$,
where $\mathcal{B}\left(\mathbb{X}\right)$ is the Borel $\sigma$-algebra
of $\mathbb{X}$. Suppose that we observe an independent and identically
distributed (IID) sequence $\left(\bm{X}_{i}\right)_{i\in\left[n\right]}$,
where $\bm{X}_{i}$ has measure $\text{Pr}_{\bm{X}}$, for each $i\in\left[n\right]=\left\{ 1,\dots,n\right\} $
($n\in\mathbb{N}$).

Let $\left(\mathscr{H}_{k}\right)_{k\in\left[m\right]}$ be a sequence
of \emph{hypotheses} that define functional spaces of \emph{models}
determined by a parameter vector $\bm{\theta}_{k}\in\mathbb{T}_{k}\subset\mathbb{R}^{q_{k}}$,
where $q_{k}\in\mathbb{N}$, for each $k\in\left[m\right]$. That
is, we can write:
\[
\mathscr{H}_{k}=\left\{ h_{k}\left(\cdot;\bm{\theta}_{k}\right):\mathbb{X}\rightarrow\mathbb{R}\text{, }\bm{\theta}_{k}\in\mathbb{T}_{k}\right\} \text{.}
\]
Further, define a \emph{loss function} $\ell:\mathbb{R}\rightarrow\mathbb{R}$
and, for each $k\in\left[m\right]$, call
\[
R_{k,n}\left(\theta_{k}\right)=\frac{1}{n}\sum_{i=1}^{n}\ell\left(h_{k}\left(\bm{X}_{i};\bm{\theta}_{k}\right)\right)
\]
and 
\[
r_{k}\left(\theta_{k}\right)=\text{E}\left\{ \ell\left(h_{k}\left(X;\theta_{k}\right)\right)\right\} 
\]
the \emph{empirical }and\emph{ expected risks}, defined by parameter
$\bm{\theta}_{k}$, within hypothesis $\mathscr{H}_{k}$, respectively.
For brevity, we will write $\ell_{k}\left(\bm{X};\bm{\theta}_{k}\right)$
in place of $\ell\left(h_{k}\left(\bm{X};\bm{\theta}_{k}\right)\right)$.

Define 
\[
k^{*}=\min\underset{k\in\left[m\right]}{\arg\min}\left\{ \min_{\theta_{k}\in\mathbb{T}_{k}}r_{k}\left(\bm{\theta}_{k}\right)\right\} 
\]
to be the \emph{optimal hypothesis index}. Here, it is assumed that
there is an order of `\emph{complexity}' that is indexed by $k$,
such as when $\left(\mathscr{H}_{k}\right)_{k\in\left[m\right]}$
is a nested sequence (i.e., $\mathscr{H}_{k}\subset\mathscr{H}_{k+1}$,
for each $k\in\left[m-1\right]$), or when the parameter spaces $\left(\mathbb{T}_{k}\right)_{k\in\left[m\right]}$
have increasing dimensions (i.e., $q_{k}<q_{k+1}$, for each $k\in\left[m-1\right]$),
which may correspond to both nested and non-nested sequences of hypotheses.
Our objective is to obtain a sequence of estimates $\left(\hat{K}_{n}\right)_{n\in\mathbb{N}}$
of $k^{*}$, where $\hat{K}_{n}:\mathbb{X}^{n}\rightarrow\left[m\right]$,
for each $n\in\mathbb{N}$, and where $\hat{K}_{n}$ has the consistency
property:
\begin{equation}
\lim_{n\rightarrow\infty}\text{Pr}\left(\hat{K}_{n}=k^{*}\right)=1\text{.}\label{eq: Consistency}
\end{equation}

The problem of estimating $k^{*}$ has long been studied in statistics
and machine learning under the topics of model identification, model
selection, variable selection, and structural risk minimisation, among
other names, where $\mathscr{H}_{k^{*}}$ is often referred to as
the class of \emph{parsimonious} models. We will generically refer
to the problem as \emph{model selection}, without loss of generality.
For treatments on the subject, we refer the interested reader to the
volumes of \citet{McQuarrie1998Regression-and-}, \citet{Burnham2002},
\citet{Massart2007}, \citet{Claeskens2008}, \citet{Konishi2008},
\citet{Ando2010Bayesian-Model-}, \citet{Gassiat2018}, and \citet{Oneto2020Model-Selection},
as well as the expositions of \citet[Ch. 6]{Vapnik1998}, \citet{Leeb2005Model-selection},
\citet{Leeb:2009aa}, \citet[Ch. 4]{Mohri2018Foundations-of-}, and
\citet[Ch. 28]{Hansen2022Econometrics}.

In this work, we consider the method of penalisation, best identified
with \emph{information criteria} (IC), that was pioneered in the works
of \citet{Akaike1974} and \citet{Schwarz1978}. That is, we propose
estimators of $k^{*}$ of the form 
\begin{equation}
\hat{K}_{n}=\min\underset{k\in\left[m\right]}{\arg\min}\text{ }\left\{ \min_{\bm{\theta}_{k}\in\mathbb{T}_{k}}R_{k,n}\left(\bm{\theta}_{k}\right)+P_{k,n}\right\} \text{,}\label{eq: general IC}
\end{equation}
where $P_{k,n}:\Omega\rightarrow\mathbb{R}_{\ge0}$ are random variables
dependent on the data and hypothesis index, which we shall refer to
as \emph{penalt}ies. 

The reader will be able to identify this framework with the Akaike
IC (AIC; \citealp{Akaike1974}) and the Bayesian IC (BIC; \citealp{Schwarz1978}),
respectively, by setting $\ell$ to be the negative logarithm, $\left(\mathscr{H}_{k}\right)_{k\in\left[m\right]}$
to be sets of probability density functions, and by taking $P_{k,n}$
to be equal to either $n^{-1}q_{k}$ or $n^{-1}q_{k}\log n/2$. 

Let $\bm{\theta}_{k}^{*}$ be a minimiser of $r_{k}$, for each hypothesis
$\mathscr{H}_{k}$. In \citet{Sin:1996aa}, general consistency results
(i.e., sufficient conditions guaranteeing (\ref{eq: Consistency}))
are provided for broad classes of losses and penalties, under strong
assumptions regarding $\bm{\theta}_{k}^{*}$, $\ell$, and $h_{k}$.
For example, $\bm{\theta}_{k}^{*}$ is required to be unique and an
interiour point of $\mathbb{T}_{k}$, for each $n$, and $\ell_{k}\left(x;\cdot\right)$
must be twice continuously differentiable, for each $x\in\mathbb{X}$,
where $\partial^{2}r\left(\bm{\theta}\right)/\partial\bm{\theta}\partial\bm{\theta}^{\top}$
is assumed to exist and be positive definite at $\bm{\theta}_{k}^{*}$.
We note that comparable general results were obtained in \citet{Baudry:2015aa}.

We generalise the results of \citet{Sin:1996aa} and \citet{Baudry:2015aa}
by providing conditions under which we may admit non-uniqueness of
$\bm{\theta}_{k}^{*}$ and lack of differentiability of $\ell_{k}\left(x;\cdot\right)$.
This includes not requiring the positive-definiteness assumption on
$\partial^{2}r\left(\bm{\theta}_{k}^{*}\right)/\partial\bm{\theta}\partial\bm{\theta}^{\top}$,
stated above. This permits simple constructions of information criteria
for problems where such assumptions are non-trivial, such as for finite
mixture models, as we discuss in the sequel. Other models where such
pathologies occur are the so-called singular models discussed in \citet{Watanabe2009}
and \citet{drton2017bayesian}, which include common models such as
factor analyzers, Bayesian networks and neural networks with hidden
nodes, reduced rank regression models, and variants of such models,
among others. This is achieved via asymptotic results for sample average
approximations of stochastic programming problems, as covered in \citet{Shapiro2000}
and \citet[Ch. 5]{Shapiro2021}. Due to the generality of our result
for consistent estimation of $k^{*}$ across different modeling and
learning problems, and due to our use of the IC form (\ref{eq: general IC}),
we name our approach \emph{PanIC}, using the Greek prefix \emph{`pan'}
to mean \emph{`all'} or \emph{`of everything'}.

Consistency results for IC have been studied in specifics in numerous
previous works. For example, one may point to the thorough analyses
of \citet{Leroux1992}, \citet{Keribin2000}, \citet{Hui:2015aa},
\citet{drton2017bayesian}, and \citet[Ch. 4]{Gassiat2018}, when
$\ell$ is taken as the negative logarithm and $\mathscr{H}_{k}$
are classes of finite mixture models; \citet{Varin2005}, \citet{Gao2010},
\citet{Ng2014a}, and \citet{Hui2021}, when $\ell_{k}$ are taken
to be the negative composite, pseudo-, or quasi-log-likelihood objects;
and \citet{Kobayashi2006}, \citet{Claeskens2008a}, and \citet{Zhang2016},
in the support vector machine context. Furthermore, the finite sample
properties of the IC has been studied in many situations, for example,
in the works of \citet{Barron1999}, \citet{Massart2007}, \citet{Bartlett2008},
and \citet{Giraud2022}. Instead of seeking to replace the thorough
and specific treatments provided by the listed texts, we aim to present
general results for obtaining consistent estimators of $k^{*}$ in
situations where other approaches may be difficult to verify or are
unavailable.

The rest of the work proceeds as follows. We provide an exposition
of the PanIC approach and its theoretical properties in Section 2.
Example applications are provided in Section 3 and numerical experiments
are performed in Section 4. We discuss our work in Section 5. Proofs
and technical results are relegated to the Appendix of the text.

\section{Main results}

We shall retain the notation and technical setting from Section \ref{sec:Introduction}
and make the following assumptions, for each $k\in\left[m\right]$:
\begin{description}
\item [{A1}] $\mathbb{T}_{k}$ is compact and there exists a $\bm{\tau}_{k}\in\mathbb{T}_{k}$,
such that
\[
\text{E}\left\{ \ell_{k}\left(\bm{X};\bm{\tau}_{k}\right)^{2}\right\} <\infty\text{.}
\]
\item [{A2}] There exists a measurable function $\mathfrak{C}_{k}:\mathbb{X}\rightarrow\mathbb{R}_{\ge0}$,
such that $\text{E}\left\{ \mathfrak{C}_{k}\left(\bm{X}\right)^{2}\right\} <\infty$
and
\[
\left|\ell_{k}\left(\bm{x};\bm{\theta}_{k}\right)-\ell_{k}\left(\bm{x};\bm{\tau}_{k}\right)\right|\le\mathfrak{C}_{k}\left(\bm{x}\right)\left\Vert \bm{\theta}_{k}-\bm{\tau}_{k}\right\Vert \text{,}
\]
for every $\bm{\theta}_{k},\bm{\tau}_{k}\in\mathbb{T}_{k}$ and almost
every $\bm{x}\in\mathbb{X}$.
\end{description}
Notice that A1 and A2 together imply that $\ell_{k}\left(\bm{x};\bm{\theta}_{k}\right)$
is \emph{Caratheodory} in the sense that $\ell_{k}\left(\bm{x};\cdot\right):\mathbb{T}_{k}\rightarrow\mathbb{R}$
is continuous for each $\bm{x}\in\mathbb{X}$, and $\ell_{k}\left(\cdot;\bm{\theta}_{k}\right):\mathbb{X}\rightarrow\mathbb{R}$
is measurable for each $\bm{\theta}_{k}\in\mathbb{T}_{k}$, for each
$k\in\left[m\right]$. We write
\[
\mathbb{T}_{k}^{*}=\arg\min_{\bm{\theta}_{k}\in\mathbb{T}_{k}}\text{ }r_{k}\left(\bm{\theta}_{k}\right)
\]
to denote the set of minima of the expected risk under hypothesis
$\mathscr{H}_{k}$. We shall indicate \emph{convergence in distribution}
by the symbol $\rightsquigarrow$. A direct application of \citet[Thm. 5.7]{Shapiro2021}
yields the following lemma, which is instrumental in the proof of
the main result in the following section.
\begin{lem}
\label{Lem: A1-A3}Assume that $\left(\bm{X}_{i}\right)_{i\in\left[n\right]}$
is an IID sequence, and that A1 and A2 hold, for each $k\in\left[m\right]$.
Then, for each $k$,
\begin{equation}
\sqrt{n}\left(\min_{\bm{\theta}_{k}\in\mathbb{T}_{k}}R_{k,n}\left(\theta_{k}\right)-\min_{\bm{\theta}_{k}\in\mathbb{T}_{k}}r_{k}\left(\bm{\theta}_{k}\right)\right)\rightsquigarrow\inf_{\bm{\theta}_{k}\in\mathbb{T}_{k}^{*}}Y_{k}\left(\cdot;\bm{\theta}_{k}\right)\text{,}\label{eq: Lemma 1}
\end{equation}
where $Y_{k}:\Omega\times\mathbb{T}_{k}\to\mathbb{R}$ is a bounded
random process, indexed by $\mathbb{T}_{k}$, with continuous sample
paths $Y_{k}\left(\omega;\cdot\right)$. 
\end{lem}
\begin{rem}
We note that details regarding $Y_{k}$ are immaterial, since $Y_{k}\left(\omega;\bm{\theta}_{k}\right)$
must be a bounded and continuous function of $\bm{\theta}_{k}\in\mathbb{T}_{k}$,
and thus $\inf_{\theta_{k}\in\mathbb{T}_{k}^{*}}Y_{k}\left(\omega;\bm{\theta}_{k}\right)$
is finite for every $\omega\in\Omega$. Therefore $\inf_{\bm{\theta}_{k}\in\mathbb{T}_{k}^{*}}Y_{k}\left(\cdot;\bm{\theta}_{k}\right):\Omega\rightarrow\mathbb{R}$
is a random variable in the usual sense. Since the left-hand side
(LHS) of (\ref{eq: Lemma 1}) converges in distribution to a random
variable, this permits us to conclude that the LHS of (\ref{eq: Lemma 1})
is bounded in probability, as $n\rightarrow\infty$ (cf. \citealt[Sec. 5.5.3c]{Boos2013}).
\end{rem}

\subsection{\label{subsec:PanIC}PanIC}

Make the following additional assumptions, for each $k\in\left[m\right]$:
\begin{description}
\item [{B1}] $P_{k,n}>0$, for each $k\in\left[m\right]$ and $n\in\mathbb{N}$,
and $P_{k,n}=o_{\text{Pr}}\left(1\right)$, as $n\rightarrow\infty$.
\item [{B2}] If $k<l$, then $\sqrt{n}\left\{ P_{l,n}-P_{k,n}\right\} \rightarrow\infty$,
in probability, as $n\rightarrow\infty$.
\end{description}
We will say that any estimator of form (\ref{eq: general IC}), satisfying
B1 and B2, is a \emph{PanIC}. The proof of the following result appears
in the Appendix, along with those of Theorems \ref{Theorem: BIC-like}
and \ref{Theorem: SWIC a.s.}.

\begin{thm}
\label{Theorem: SWIC}Assume that $\left(\bm{X}_{i}\right)_{i\in\left[n\right]}$
is an IID sequence, and that A1, A2, B1 and B2 hold, for each $k\in\left[m\right]$.
Then, the PanIC (\ref{eq: general IC}) satisfies the consistency
property (\ref{eq: Consistency}).
\end{thm}
\begin{rem}
\label{Rem: SWIC comparison}We note that PanIC were considered in
\citet[Prop. 4.6(a,b)]{Sin:1996aa}. However, as noted in Section
\ref{sec:Introduction}, the \citet{Sin:1996aa} results required
much stronger assumptions than A1 and A2. Namely, it is assumed that
both $\ell_{k}\left(\bm{x};\cdot\right)$ and $r_{k}$ are both twice
continuously differentiable on a set $\mathbb{T}_{k}$, with non-empty
interiour, for each $k\in\left[m\right]$; that $r_{k}$ has a unique
minimiser at $\bm{\theta}_{k}^{*}$, in the interiour of $\mathbb{T}_{k}$,
for each $n\in\mathbb{N}$; that $\ell_{k}\left(\bm{X};\cdot\right)$
and its derivative $\partial\ell_{k}\left(\bm{X};\bm{\theta}_{k}\right)/\partial\bm{\theta}_{k}$
both satisfy central limit theorems; and the Hessian $\partial^{2}r_{k}\left(\bm{\theta}_{k}\right)/\partial\bm{\theta}_{k}\partial\bm{\theta}_{k}^{\top}$
is positive definite at $\bm{\theta}_{k}^{*}$. In contrast, Theorem
\ref{Theorem: SWIC} makes no assumptions regarding the uniqueness
or even countability of minima of $r_{k}$ and $R_{k,n}$, nor the
location of the minima within $\mathbb{T}_{k}$. Moreover, the totality
of A1 and A2 amount to a Lipschitz assumption on $\ell_{k}\left(\bm{x};\cdot\right)$,
and thus demands less than one derivative on each of the $k$ functions,
although if $\ell_{k}\left(\bm{x};\cdot\right)$ are each continuously
differentiable then the verification of A1 and A2 becomes procedural,
as one can take the supremum norm of the gradients, with respect to
$\bm{\theta}_{k}$, in place of $\mathfrak{C}_{k}\left(\bm{x}\right)$. 
\end{rem}
\begin{rem}
\label{rem: SWIC}Let $\left(c_{k}\right)_{k\in\left[m\right]}$ be
positive constants, such that $c_{k}<c_{l}$, for each $k<l$, and
define $\log_{+}\left(\cdot\right)=\max\left\{ 1,\log\left(\cdot\right)\right\} $.
We are particularly interested in penalties taking the form:
\begin{equation}
P_{n,k}=\alpha c_{k}\sqrt{\frac{\log_{+}^{\left(\beta\right)}\left(n\right)}{n}}\text{,}\label{eq: SWIC pen}
\end{equation}
where $\alpha>0$, $\beta\in\mathbb{N}$, and
\[
\log_{+}^{\left(\beta\right)}\left(n\right)=\underset{\beta\text{ times}}{\underbrace{\left(\log_{+}\circ\log_{+}\circ\dots\circ\log_{+}\right)}}\left(n\right)\text{.}
\]
We call any IC with penalty of form (\ref{eq: SWIC pen}) a \emph{SWIC}
(Sin--White IC) of order $\beta$, with constant $\alpha$, in homage
to \citet{Sin:1996aa}, who suggest the use of the first and second
order versions.

We can easily check that penalties of form (\ref{eq: SWIC pen}) satisfy
B1 and B2, and thus SWIC are within the class of PanIC. B1 is verified
by noting that $\log_{+}^{\left(\beta\right)}\left(n\right)$ is increasing
and real-valued and increases slower than $\sqrt{n}$. B2 can be checked
by observing that
\[
\sqrt{n}\left\{ P_{l,n}-P_{k,n}\right\} =\alpha\left(c_{l}-c_{k}\right)\sqrt{\log_{+}^{\left(\beta\right)}\left(n\right)}\underset{n\rightarrow\infty}{\longrightarrow}\infty\text{,}
\]
since $\alpha\left(c_{l}-c_{k}\right)>0$, for $k<l$, as required.
\end{rem}

\subsection{BIC-like criteria}

We shall say that a criterion is \emph{BIC-like} if it satisfies the
following assumption, along with B1, for each $k\in\left[m\right]$:
\begin{description}
\item [{B2$^{*}$}] If $k<l$, then $n\left\{ P_{l,n}-P_{k,n}\right\} \rightarrow\infty$,
in probability, as $n\rightarrow\infty$.
\end{description}
General conditions under which BIC-like estimators (\ref{eq: general IC})
are consistent were also studied in \citet[Prop. 4.2(c)]{Sin:1996aa}
and \citet[Thm. 8.1 and Cor. 8.2]{Baudry:2015aa}. Using further results
from the stochastic programming literature, we can obtain alternative
consistency criteria that are not implied by previous works. We make
the following additional assumptions, for each $k\in\left[m\right]$:
\begin{description}
\item [{C1}] $r_{k}$ is Lipschitz continuous on $\mathbb{T}_{k}$, and
is uniquely minimised and twice differentiable at some $\bm{\theta}_{k}^{*}\in\mathbb{T}_{k}$.
\item [{C2}] $\ell_{k}\left(\bm{x};\cdot\right)$ is Lipschitz continuous
and differentiable at $\bm{\theta}_{k}^{*}$, for almost all $x\in\mathbb{X}$.
\item [{C3}] The set $\mathbb{T}_{k}$ is \emph{second order regular} at
$\bm{\theta}_{k}^{*}$ (see the Appendix).
\item [{C4}] The \emph{quadratic growth condition} holds at $\bm{\theta}_{k}^{*}$
(see the Appendix).
\item [{C5}] If $r\left(\theta_{k}^{*}\right)=r\left(\theta_{k^{*}}^{*}\right)$,
then $n\left(R_{k,n}\left(\theta_{k}^{*}\right)-R_{k^{*},n}\left(\theta_{k^{*}}^{*}\right)\right)=O_{\text{Pr}}\left(1\right)$.
\end{description}
\begin{thm}
\label{Theorem: BIC-like}Assume that $\left(\bm{X}_{i}\right)_{i\in\left[n\right]}$
is an IID sequence, and that A1, A2, B1, B2$^{*}$, and C1--C5 hold,
for each $k\in\left[m\right]$. Then, the BIC-like estimator (\ref{eq: general IC})
satisfies the consistency property (\ref{eq: Consistency}).
\end{thm}
\begin{rem}
\citet[Prop. 4.2(c)]{Sin:1996aa} prove the consistency of BIC-like
criteria under the same assumptions as those discussed in Remark \ref{Rem: SWIC comparison}.
A slight relaxation is provided in \citet[Thm. 8.1 and Cor. 8.2]{Baudry:2015aa},
where each $\ell_{k}\left(\bm{x};\cdot\right)$ is only assumed to
be continuously differentiable on $\mathbb{T}_{k}$, for almost all
$\bm{x}\in\mathbb{X}$, for $k\in\left[m\right]$. However, it is
still assumed that $\bm{\theta}_{k}^{*}$ is in the interiour of $\mathbb{T}_{k}$,
and that the Hessian $\partial^{2}r_{k}\left(\bm{\theta}_{k}\right)/\partial\bm{\theta}_{k}\partial\bm{\theta}_{k}^{\top}$
is positive definite at $\bm{\theta}_{k}^{*}$, for each $k$. From
Appendix A, we observe that Theorem \ref{Theorem: BIC-like} makes
the same conclusions as the result of \citet{Baudry:2015aa}, under
the same assumptions on $\ell_{k}$ and $r_{k}$. However Theorem
\ref{Theorem: BIC-like} provides further relaxations, allowing $\ell\left(\bm{x};\cdot\right)$
to be Lipschitz on $\mathbb{T}_{k}\backslash\left\{ \bm{\theta}_{k}^{*}\right\} $,
and only requiring differentiability at $\bm{\theta}_{k}^{*}$, for
almost all $\bm{x}\in\mathbb{X}$. Furthermore, $\bm{\theta}_{k}^{*}$
is permitted to be on the boundary of $\mathbb{T}_{k}$, with C3 and
C4 being particularly easy to verify when $\mathbb{T}_{k}$ are polyhedral
sets. We note that the result of \citet{Sin:1996aa} is proved by
first verifying that $\bm{\theta}_{k}^{*}$ is asymptotically normal
for each $k\in\left[m\right]$. An asymptotic expansion is then used
to obtain the boundedness in probability of the $n$-scaled minimum
risk. Alternatively, \citet{Baudry:2015aa} takes an empirical processes
approach, where concentration inequalities for relevant terms are
obtained and used to obtain the boundedness in probability of required
terms.
\end{rem}
\begin{rem}
In \citet[Cor. 8.2]{Baudry:2015aa}, it is proposed that C5 can be
implied by the more explicit assumption that $r_{k}\left(\bm{\theta}_{k}^{*}\right)=r_{k^{*}}\left(\bm{\theta}_{k^{*}}^{*}\right)$
if and only if $\ell_{k}\left(\bm{x};\bm{\theta}_{k}\right)=\ell_{k^{*}}\left(\bm{x};\bm{\theta}_{k^{*}}^{*}\right)$,
for almost all $x\in\mathbb{X}$. \citet{Sin:1996aa} do not propose
any method for checking C5 but suggests that it is generally satisfied
in situations where the hypotheses are nested, in the sense that $\mathscr{H}_{k}\subset\mathscr{H}_{l}$,
for each $k<l$, particularly when $\ell$ is taken to be the negative
log-likelihood and $\mathscr{H}_{k}$ are sets of probability density
functions, for each $k\in\left[m\right]$ (cf. \citealp{Vuong1989Likelihood-rati}). 
\end{rem}
\begin{rem}
Of course our reference to ICs with penalties satisfying B1 and B2$^{*}$
as BIC-like is suggestive that the BIC should fall within this category.
Upon assuming that $q_{k}<q_{l}$ for each $k<l$, the original BIC
penalty of \citet{Schwarz1978} corresponds to the choice $P_{k,n}=q_{k}n^{-1}\log n/2$,
which verifies B2$^{*}$ by observing that
\[
n\left\{ P_{l,n}-P_{k,n}\right\} =\left(q_{l}-q_{k}\right)\log n/2\underset{n\rightarrow\infty}{\longrightarrow}\infty\text{,}
\]
since $q_{l}-q_{k}>0$. Note that the class of BIC-like criteria is
not limited to the BIC and include, for example, the Hannan--Quinn
information criterion with penalty $P_{k,n}=q_{k}n^{-1}\log_{+}^{\left(2\right)}\left(n\right)/2$
\citep{Hannan1979}, or more generally, IC with penalties of the form
$P_{k,n}=\alpha q_{k}n^{-1}\log_{+}^{\left(\beta\right)}\left(n\right)$,
for $\alpha>0$ and $\beta\in\mathbb{N}$.
\end{rem}

\section{Example applications}

\subsection{Finite mixture models}

For each $k\in\left[m\right]$, let
\begin{equation}
\mathscr{H}_{k}=\left\{ \begin{array}{c}
h_{k}\left(\cdot;\bm{\theta}_{k}\right)=\sum_{z=1}^{k}\pi_{z}f\left(\cdot;\bm{\upsilon}_{z}\right):\theta_{k}=\left(\pi_{1},\dots,\pi_{k},\bm{\upsilon}_{1},\dots,\bm{\upsilon}_{k}\right)\text{,}\\
\theta_{k}\in\mathbb{T}_{k}=\mathbb{S}_{k-1}\times\mathbb{U}^{k}
\end{array}\right\} \text{,}\label{eq: mixture hypotheses}
\end{equation}
where $f\left(\cdot;\bm{\upsilon}\right)$ is a (\emph{component})
density with parameter $\bm{\upsilon}\in\mathbb{U}\subset\mathbb{R}^{q}$,
where $\mathbb{U}$ is compact and $\mathbb{S}_{k-1}$ is the probability
simplex in $\mathbb{R}^{k}$. We call $\mathscr{H}_{k}$ the set of
\emph{finite mixtures} of component densities $f\left(\bm{x};\bm{\upsilon}\right)$
of order $k$.

Suppose that we observe IID data $\left(\bm{X}_{i}\right)_{i\in\left[n\right]}$
from the density
\[
h_{k^{*}}\left(\bm{x};\bm{\theta}_{k^{*}}^{*}\right)=\sum_{z=1}^{k^{*}}\pi_{z}^{*}f\left(\bm{x};\bm{\upsilon}_{z}^{*}\right)\in\mathscr{H}_{k^{*}}\text{,}
\]
for some $k^{*}\in\left[m\right]$, where $\bm{\theta}_{k}^{*}=\left(\pi_{1}^{*},\dots,\pi_{k}^{*},\bm{\upsilon}_{1}^{*},\dots,\bm{\upsilon}_{k}^{*}\right)\in\mathbb{T}_{k}$,
for each $k\in\left[m\right]$, and
\begin{equation}
\bm{\theta}_{k}^{*}\in\underset{\bm{\theta}_{k}\in\mathbb{T}_{k}}{\arg\min}\text{ }\underset{=r_{k}\left(\bm{\theta}_{k}\right)}{\underbrace{\text{E}\left\{ -\log h_{k}\left(\bm{X};\bm{\theta}_{k}\right)\right\} }}\text{.}\label{eq: finite mixture optim}
\end{equation}
In the classical IC setting, if $k^{*}$ is unknown, then we typically
estimate it via a penalised estimator of the form
\[
\hat{K}_{n}=\min\underset{k\in\left[m\right]}{\arg\min}\text{ }\left\{ \min_{\theta_{k}\in\mathbb{T}_{k}}\underset{=R_{k,n}\left(\bm{\theta}_{k}\right)}{\underbrace{-\frac{1}{n}\sum_{i=1}^{n}\log h_{k}\left(\bm{X}_{i};\bm{\theta}_{k}\right)}}+P_{k,n}\right\} \text{,}
\]
where $R_{k,n}\left(\bm{\theta}_{k}\right)$ is the average negative
log-likelihood the model in $\mathscr{H}_{k}$ defined by parameter
$\bm{\theta}_{k}\in\mathbb{T}_{k}$. This is often referred to as
the problem of \emph{order selection} in finite mixture modeling. 

To fulfill B1 and B2, we can propose SWIC-type penalities of the form
\begin{align}
P_{k,n} & =\alpha q_{k}\sqrt{n^{-1}\log^{\left(\beta\right)}\left(n\right)}=\alpha k\left(1+q\right)\sqrt{n^{-1}\log^{\left(\beta\right)}\left(n\right)}\text{,}\label{eq: SWIC mixture}
\end{align}
and in the case of B2$^{*}$, we can propose the typical BIC penalty
\begin{align}
P_{k,n} & =\frac{q_{k}}{2}\frac{\log n}{n}=\frac{k\left(1+q\right)}{2}\frac{\log n}{n}\text{.}\label{eq: BIC mixture}
\end{align}

Applying Theorem \ref{Theorem: SWIC} is typically procedural, as
A1 and A2 can also be checked directly. For example, supposing that
we have a mixture of normal component densities:
\[
h_{k}\left(x;\bm{\theta}_{k}\right)=\sum_{z=1}^{k}\pi_{k}\phi\left(x;\mu_{k},\sigma_{k}^{2}\right)\text{,}
\]
where 
\[
\phi\left(x;\mu,\sigma^{2}\right)=\frac{1}{\sqrt{2\pi\sigma^{2}}}\exp\left\{ -\frac{1}{2}\frac{\left(x-\mu\right)^{2}}{\sigma^{2}}\right\} 
\]
is the normal density with mean $\mu\in\mathbb{R}$ and variance $\sigma^{2}>0$,
$\bm{\theta}_{k}=\left(\pi_{1},\dots,\pi_{k},\mu_{1},\sigma_{1}^{2},\dots,\mu_{k},\sigma_{k}^{2}\right)$,
and $\mu_{k}\in\left[-\mathfrak{m},\mathfrak{m}\right]$ and $\sigma_{k}^{2}\in\left[1/\mathfrak{s},\mathfrak{s}\right]$
for constants $\mathfrak{m}>0$ and $\mathfrak{s}>1$. Then, we have
the fact that $h\left(\bm{x};\cdot\right)$ is differentiable, for
each $x\in\mathbb{R}$, with gradient:
\begin{align*}
 & \frac{\partial\log h_{k}\left(x;\bm{\theta}_{k}\right)}{\partial\bm{\theta}_{k}}\\
= & \left(\dots,\frac{\phi\left(x;\mu_{k},\sigma_{k}^{2}\right)}{h\left(x;\bm{\theta}_{k}\right)},\frac{x-\mu_{1}}{\sigma_{1}^{2}}\pi_{1}\frac{\phi\left(x;\mu_{1},\sigma_{1}^{2}\right)}{h\left(x;\bm{\theta}_{k}\right)},\frac{\left(\mu_{1}-x\right)^{2}-\sigma_{1}^{2}}{2\sigma_{1}^{4}}\pi_{1}\frac{\phi\left(x;\mu_{1},\sigma_{1}^{2}\right)}{h\left(x;\bm{\theta}_{k}\right)},\dots\right)
\end{align*}
and so $-\log h_{k}\left(x;\cdot\right)$ is Lipschitz continuous
with a constant of the form 
\[
\mathfrak{C}\left(x\right)=k\times\max_{\mu\in\left[-\mathfrak{m},\mathfrak{m}\right],\sigma^{2}\in\left[1/\mathfrak{s},\mathfrak{s}\right]}\max\left\{ 1,\frac{\left|x-\mu\right|}{\sigma^{2}},\frac{\left(x-\mu\right)^{2}+\sigma^{2}}{2\sigma^{4}}\right\} \text{,}
\]
since $\phi\left(x;\mu_{z},\sigma_{z}^{2}\right)/h_{k}\left(x;\bm{\theta}_{k}\right)$
and $\pi_{k}$ are bounded by 1, for each $z\in\left[k\right]$, via
the mean value theorem. Then, we can verify that A3 holds by noting
that normal random variables have finite moments of every order. A2
can be verified by the same observation.

For general sets $\mathbb{T}_{k}$, we cannot verify C1 since finite
mixture models are only identifiable up to a permutation of the parameters.
That is, if $-\log h\left(x;\cdot\right)$ is minimised at $\bm{\theta}_{k}^{*}$,
then it is also minimised at
\[
\bm{\theta}_{k}^{\Pi}=\left(\pi_{\Pi\left(1\right)}^{*},\dots,\pi_{\Pi\left(k\right)}^{*},\mu_{\Pi\left(1\right)}^{*},\sigma_{\Pi\left(1\right)}^{*2},\dots,\mu_{\Pi\left(k\right)}^{*},\sigma_{\Pi\left(k\right)}^{*2}\right)\text{,}
\]
where $\Pi:\left[k\right]\rightarrow\left[k\right]$ is a permutation.
However, finite mixture models can often have isolated minima (\ref{eq: finite mixture optim})
and thus we can always consider $\mathbb{T}_{k}$ to be a compact
set that contains only one such optimal parameter. C2 is also more
onerous than A3, but is can be verified in some situations. C3 and
C4 can be difficult to verify, given that mixture models often have
singular Fisher information (see, e.g. \citealp[Sec. 7.8]{Watanabe2009}).
Thus, it is difficult to make use of Theorem \ref{Theorem: BIC-like}
for finite mixture models. However, this does not dismiss the use
of BIC-like IC, as demonstrated by the specific analyses of \citet{Leroux1992},
\citet{Keribin2000}, and \citet[Ch. 4]{Gassiat2018}, who establish
conditions under which BIC-like IC can be proved consistent for mixture
model order selection.
\begin{rem}
\label{Rem: Parametric vs nonparametric}We note that since $\left(\mathscr{H}_{k}\right)_{k\in\left[m\right]}$
is a nested sequence of hypotheses, our problem is only nontrivial
if we assume that for sufficiently large $m$, there exists some $k^{*}<m$
such that $\inf_{\bm{\theta}_{k^{*}}\in\mathbb{T}_{k^{*}}}r_{k^{*}}\left(\bm{\theta}_{k^{*}}\right)=\inf_{\bm{\theta}_{k}\in\mathbb{T}_{k}}r_{k}\left(\bm{\theta}_{k}\right)$,
for any $k>k^{*}$. This is true, for example, in the well-specified
case, where the underlying data generating process is has density
function in $\mathscr{H}_{k}$, for every $k\ge k^{*}$. That is,
the data $\left(\bm{X}_{i}\right)_{i\in\left[n\right]}$ are sampled
IID from a normal mixture model with $k^{*}\in\left[m\right]$ components.
If this is not the case, in the sense that $\inf_{\bm{\theta}_{k}\in\mathbb{T}_{k}}r_{k}\left(\bm{\theta}_{k}\right)$
is a strictly decreasing function of $k\in\mathbb{N}$, then the order
selection problem above no longer makes sense and we then seek instead
to choose a map $n\mapsto k\left(n\right)$, such that the loss expected
loss $\text{E}\left[\inf_{\bm{\theta}_{k\left(n\right)}\in\mathbb{T}_{k\left(n\right)}}R_{k\left(n\right),n}\left(\bm{\theta}_{k\left(n\right)}\right)\right]$
converges to the limit $\lim_{k\rightarrow\infty}\inf_{\bm{\theta}_{k}\in\mathbb{T}_{k}}r_{k}\left(\bm{\theta}_{k}\right)$
at the fastest rate. For the minimum average negative log-likelihood
problem, above, such results are considered in \citet{LiBarron1999}
and \citet{RakhlinPanchenkoMukherjee2005}. In particular, when the
assumptions of \citet{RakhlinPanchenkoMukherjee2005} hold, the optimal
rate of should take $k\left(n\right)=O\left(\sqrt{n}\right)$, which
yields a convergence rate of $O\left(1/\sqrt{n}\right)$. Examples
of recent and more sophisticated results in this direction can be
found in \citet{ho2022convergence}.
\end{rem}

\subsection{Least absolute deviation and $\epsilon$-support vector regression\label{subsec:Least-absolute-deviation}}

We suppose that we observe IID data $\left(\bm{X}_{i}\right)_{i\in\left[n\right]}$,
where each $\bm{X}_{i}$ has the same distribution as $\bm{X}=\left(Y,\bm{W}\right)\in\mathbb{X}=\mathbb{Y}\times\mathbb{W}$,
$\mathbb{Y}\subset\mathbb{R}$, and $\mathbb{W}\subset\mathbb{R}^{m}$.
For each $k\in\left[m\right]$, we let
\begin{equation}
\mathscr{H}_{k}=\left\{ h_{k}\left(\bm{x};\bm{\theta}_{k}\right)=y-\bm{w}^{\left(k\right)\top}\bm{\theta}_{k}:\bm{\theta}_{k}\in\mathbb{T}_{k}\right\} \label{eq: SVR hypo}
\end{equation}
be the \emph{residual function }obtained from the estimation of $y\in\mathbb{R}$
by $\bm{w}^{\left(k\right)\top}\bm{\theta}_{k}$, where $\left(\cdot\right)^{\top}$
is the transposition operator, and $\bm{w}^{\left(k\right)}=\left(w_{1},\dots,w_{k}\right)\in\mathbb{R}^{k}$
is the first $k$ elements of the vector $\bm{w}\in\mathbb{R}^{m}$.
Here, we assume that $\mathbb{T}_{k}\subset\mathbb{R}^{k}$ is compact,
for each $k\in\left[m\right]$. We shall measure the goodness of the
approximation $\bm{w}^{\left(k\right)\top}\bm{\theta}_{k}$, for each
$k$, using the \emph{$\epsilon$-insensitive $\text{L}_{1}$ loss
}
\[
\ell_{\epsilon}\left(y\right)=\left(\left|y\right|-\epsilon\right)_{+}\text{,}
\]
where $\epsilon\ge0$ and $\left(\cdot\right)_{+}=\max\left\{ 0,\cdot\right\} $.
When $\epsilon=0$, we simply have the $\text{L}_{1}$ loss: $\ell_{0}\left(y\right)=\left|y\right|$.
For hypothesis class $\mathscr{H}_{k}$, we define
\[
\mathbb{T}_{k}^{*}=\underset{\bm{\theta}_{k}\in\mathbb{T}}{\arg\min}\text{ }\underset{=r_{k}\left(\bm{\theta}_{k}\right)}{\underbrace{\text{E}\left\{ \ell_{\epsilon}\left(h_{k}\left(\bm{X};\bm{\theta}_{k}\right)\right)\right\} }}
\]
to be the set of optimal linear \emph{$\epsilon$-support vector regression}
($\epsilon$-SVR) parameters, with typical elements $\bm{\theta}_{k}^{*}$,
which we can estimate empirically by an $\epsilon$-SVR estimator
(cf. \citealp[Sec. 9.1]{ScholkopfSmola2002}):
\begin{equation}
\hat{\bm{\theta}}_{k,n}\in\underset{\bm{\theta}_{k}\in\mathbb{T}}{\arg\min}\text{ }\frac{1}{n}\sum_{i=1}^{n}\ell_{\epsilon}\left(h_{k}\left(\bm{X}_{i};\bm{\theta}_{k}\right)\right)\text{.}\label{eq: SVR esti}
\end{equation}
When $\epsilon=0$, (\ref{eq: SVR esti}) is commonly referred to
as a \emph{least absolute deviation} (LAD) estimator. Supposing that
\[
k^{*}=\min\underset{k\in\left[m\right]}{\arg\min}\left\{ \min_{\bm{\theta}_{k}\in\mathbb{T}_{k}}\text{E}\left\{ \ell_{\epsilon}\left(h_{k}\left(\bm{X};\bm{\theta}_{k}\right)\right)\right\} \right\} \text{,}
\]
which we wish to estimate. This is a version of the classical \emph{variable
selection} problem in regression analysis. The typical best-subset
type variable selection problem (see, e.g., \citealp{heinze2018variable})
can be obtained by instead considering $\mathscr{H}_{k}$ to include
all models that depend on $k$ entries of $w$, rather than only the
model that depends on the first $k$ entries of $w$ as per (\ref{eq: SVR hypo}).
One can then use the suggest IC to choose the optimal number of parameters
$k^{*}$, and then within the models in $\mathscr{H}_{k^{*}}$, where
all models are equally parsimonious, one can compare the respective
sample minimum risks in order to decide which of the models with $k^{*}$
covariates is most optimal.

We can suggest IC for estimating $k^{*}$ of form
\begin{equation}
\hat{K}_{n}=\min\underset{k\in\left[m\right]}{\arg\min}\text{ }\left\{ \min_{\bm{\theta}_{k}\in\mathbb{T}_{k}}\underset{=R_{k,n}\left(\bm{\theta}_{k}\right)}{\underbrace{-\frac{1}{n}\sum_{i=1}^{n}\ell_{\epsilon}\left(h_{k}\left(\bm{X}_{i};\bm{\theta}_{k}\right)\right)}}+P_{k,n}\right\} \text{,}\label{eq: SVR ICs}
\end{equation}
where 
\begin{equation}
P_{k,n}=\alpha k\sqrt{n^{-1}\log_{+}^{\left(\beta\right)}\left(n\right)}\text{ and }P_{k,n}=\alpha kn^{-1}\log_{+}^{\left(\beta\right)}\left(n\right)\text{,}\label{eq: SWIC/BIC}
\end{equation}
fulfill the requirements of B2 and B2$^{*}$, respectively. 

Next, we observe that $\ell_{k}\left(\bm{x};\bm{\theta}_{k}\right)=\ell_{\epsilon}\left(h_{k}\left(\bm{x};\bm{\theta}_{k}\right)\right)$
is piecewise continuously differentiable in $\bm{\theta}_{k}\in\mathbb{R}^{k}$,
for each fixed $\bm{x}\in\mathbb{R}^{k}$. Furthermore, this implies
that it is therefore Lipschitz continuous, with constant $\mathfrak{C}\left(x\right)=\mathfrak{C}\times\left\Vert \bm{w}^{\left(k\right)}\right\Vert _{1}$,
for some $\mathfrak{C}>0$ (cf. \citealp[Cor. 4.1.1]{scholtes2012introduction}).
Thus, we only require that $\text{E}\left\Vert \bm{W}\right\Vert _{2}^{2}$
is finite to verify A2. To further verify A1, we only additionally
need that $\text{E}\left\{ Y^{2}\right\} $ be finite.

For arbitrary measures on $\bm{X}$ and sets $\mathbb{T}_{k}$, we
cannot guarantee that each $r_{k}\left(\bm{\theta}_{k}\right)$ has
a unique minimiser. In fact, although trivial, one can construct examples
whereupon $\mathbb{T}_{k}$ is uncountable and $\mathbb{T}_{k}^{*}=\mathbb{T}_{k}$.
As an example, for $k=1$, suppose that $W^{\left(1\right)}=1$. Then,
we can write
\[
r_{1}\left(\theta_{1}\right)=\text{E}\left\{ \left(\left|Y-\theta_{1}\right|-\epsilon\right)_{+}\right\} \text{.}
\]
Now, suppose that $Y$ is supported on the bounded set $\mathbb{Y}=\left[-2,2\right]\subset\mathbb{R}$,
$\epsilon=4$, and $\theta_{1}\in\left[-10,10\right]=\mathbb{T}_{1}$.
Then, on $\mathbb{T}_{1}$, if $Y$ is uniformly distributed on $\mathbb{Y}$,
we have
\[
r_{1}\left(\theta_{1}\right)=\begin{cases}
-4-\theta_{1} & \text{if }-10\le\theta_{1}<-6\text{,}\\
\left(\theta_{1}+2\right)^{2}/8 & \text{if }-6\le\theta_{1}<-2\text{,}\\
0 & \text{if }-2\le\theta_{1}<2\text{,}\\
\left(\theta_{1}-2\right)^{2}/8 & \text{if }2\le\theta_{1}<6\text{,}\\
\theta_{1}-4 & \text{if }6\le\theta_{1}\le10\text{,}
\end{cases}
\]
implying that $r_{1}$ is minimized by any $\theta_{1}^{*}\in\left[-2,2\right]$.
Thus, although somewhat pathological, it is possible to invalidate
the uniqueness assumption in C1. Furthermore, we have also established
that it is possible for the minimisers of $r_{k}\left(\theta_{k}\right)$
to be connected, thus unavailing the solution of localizing $\mathbb{T}_{k}$
around an appropriate isolated minimiser. 

Of course, the construction above is an exclusive pathology to the
$\epsilon>0$ case. Unfortunately, it is still not generally possible
to guarantee uniqueness of the minimiser when $\epsilon=0$, without
making additional assumptions. For example, in \citet{Pollard1991},
it is assumed that
\[
Y=\bm{\theta}_{k}^{\top}\bm{W}^{\left(k\right)}+E\text{,}
\]
where $E$ is a random variable with median zero, and a continuous
and positive density function in a neighborhood of zero, say, to guarantee
that C1 holds for $k\in\left[m\right]$. A similar approach is taken
in \citet{knight1998limiting}, who assume differentiability at zero
of the distribution function of the noise variable $E$. The approaches
of \citet{Pollard1991} and \citet{knight1998limiting} make strong
assumptions regarding the data generating process, but both provide
asymptotic normality of the risk minimiser $\bm{\theta}_{k}^{*}$,
if it is unique. Once obtained, a delta method or asymptotic expansion
approach (such as in \citealp{Sin:1996aa}) can be used to obtain
the boundedness in probability of the $\sqrt{n}$ or $n$ scaling
of the minimum risk. This differs to our approach, which seeks to
obtain the asymptotic law for the $\sqrt{n}$ scaling of the minimum
risk, directly.

We notice that if $X$ has measure that is absolutely continuous with
respect to the Lebesgue measure on $\mathbb{X}$, then $\ell_{k}\left(\bm{X};\bm{\theta}_{k}\right)=\left(\left|Y-\bm{\theta}_{k}^{\top}\bm{W}^{\left(k\right)}\right|-\epsilon\right)_{+}$
is twice differentiable with respect to $\bm{\theta}_{k}\in\mathbb{T}_{k}$,
for almost all $\bm{X}\in\mathbb{X}$, for each $k\in\left[m\right]$
but that the second derivative will be zero for every $\bm{\theta}_{k}$.
This invalidates the necessary and sufficient condition for C3 under
C4, as per the Appendix, and thus we cannot verify the consistency
of BIC-like criteria using Theorem \ref{Theorem: BIC-like}. Nevertheless,
at least in the $\epsilon=0$ case, BIC-like criteria have been proved
consistent for variable selection under various assumptions (see,
e.g., \citealp{Bai1998Estimation-of-m} and \citealp{Lee2014Model-selection}). 

Finally, we note that our construction above is a simplified version
of the variable selection problem. The typical best-subset type variable
selection problem (see, e.g., \citealp{heinze2018variable}) can be
obtained by instead considering $\mathscr{H}_{k}$ to include all
models that depend on $k$ entries of $\bm{w}$, rather than only
the model that depends on the first $k$ entries of $\bm{w}$ as per
(\ref{eq: SVR hypo}). One can use the suggest IC to choose the optimal
number of parameters $k^{*}$, and then within the models in $\mathscr{H}_{k^{*}}$,
where all models are equally parsimonious, one can compare the respective
sample minimum risks in order to decide which of the models with $k^{*}$
covariates is most optimal.

\subsection{Principal component analysis}

We next consider $\left(\bm{X}_{i}\right)_{i\in\left[n\right]}$ IID
with each $\bm{X}_{i}$ from having the same distribution as the random
variable $\bm{X}\in\mathbb{X}\subset\mathbb{R}^{m}$. For each $k\in\left[m\right]$,
we let
\begin{equation}
\mathscr{H}_{k}=\left\{ h_{k}\left(\bm{x};\bm{\theta}_{k}\right)=\bm{x}-\bm{\theta}_{k}\bm{\theta}_{k}^{\top}\bm{x}:\bm{\theta}_{k}\in\mathbb{T}_{k}\right\} \text{,}\label{eq: PCA hypo}
\end{equation}
where $\mathbb{T}_{k}\subset\mathbb{R}^{m\times k}$. Like in Section
\ref{subsec:Least-absolute-deviation}, we can consider $h_{k}\left(\bm{x};\bm{\theta}_{k}\right)$
as the residual function under a $k\le m$ dimensional linear subspace
encoding and decoding operations defined via $\bm{\theta}_{k}^{\top}$
and $\bm{\theta}_{k}$, respectively. If we choose the quadratic loss
$\ell\left(\cdot\right)=\left\Vert \cdot\right\Vert _{2}^{2}$, then
\begin{equation}
\underset{\bm{\theta}_{k}\in\mathbb{T}_{k}}{\arg\min}\underset{=R_{k,n}\left(\bm{\theta}_{k}\right)}{\underbrace{\text{ }\frac{1}{n}\sum_{i=1}^{n}\left\Vert \bm{X}_{i}-\bm{\theta}_{k}\bm{\theta}_{k}^{\top}\bm{X}_{i}\right\Vert _{2}^{2}}}\label{eq: PCA}
\end{equation}
is the classical problem of computing the first $k$ \emph{principal
components} of the data $\left(\bm{X}_{i}\right)_{i\in\left[n\right]}$,
and obtaining solutions to (\ref{eq: PCA}) is conventionally referred
to as \emph{principal component analysis} (PCA). We can conveniently
choose any compact parameter space $\mathbb{T}_{k}$ containing $\mathbb{O}_{k}=\left\{ \bm{\theta}_{k}\in\mathbb{R}^{m\times k}:\bm{\theta}_{k}^{\top}\bm{\theta}_{k}=\mathbf{I}_{k}\right\} $,
where $\mathbf{I}_{k}$ is the identity of $\mathbb{R}^{k\times k}$,
since 
\[
\underset{\bm{\theta}_{k}\in\mathbb{R}^{m\times k}}{\arg\min}\text{ }\frac{1}{n}\sum_{i=1}^{n}\left\Vert \bm{X}_{i}-\bm{\theta}_{k}\bm{\theta}_{k}^{\top}\bm{X}_{i}\right\Vert _{2}^{2}=\underset{\bm{\theta}_{k}\in\mathbb{O}_{k}}{\arg\min}\text{ }\frac{1}{n}\sum_{i=1}^{n}\left\Vert \bm{X}_{i}-\bm{\theta}_{k}\bm{\theta}_{k}^{\top}\bm{X}_{i}\right\Vert _{2}^{2}
\]
via \citet[Lem.  23.1]{Shalev-Shwartz2014}. 

We can write
\begin{equation}
\ell_{k}\left(\bm{x};\bm{\theta}_{k}\right)=\left\Vert \bm{x}-\bm{\theta}_{k}\bm{\theta}_{k}^{\top}\bm{x}\right\Vert _{2}^{2}=\left\Vert \bm{x}\right\Vert _{2}^{2}-\text{trace}\left\{ \bm{\theta}_{k}^{\top}\bm{x}\bm{x}^{\top}\bm{\theta}_{k}\right\} \text{,}\label{eq: trace version}
\end{equation}
which yields the expression
\begin{equation}
r_{k}\left(\bm{\theta}_{k}\right)=\text{trace}\left\{ \text{E}\left(\bm{X}\bm{X}^{\top}\right)\right\} -\text{trace}\left\{ \bm{\theta}_{k}^{\top}\text{E}\left(\bm{X}\bm{X}^{\top}\right)\bm{\theta}_{k}\right\} \text{,}\label{eq: r_k for PCA}
\end{equation}
via the linearity of the trace. Clearly, when $k=m$, we can choose
$\bm{\theta}_{k}^{*}=\text{\textbf{I}}_{m}$ to obtain the minimum
value $r_{k}\left(\bm{\theta}_{k}^{*}\right)=0$. However, suppose
that we wish to search for the smallest $k^{*}<m$ such that $r_{k^{*}}\left(\bm{\theta}_{k^{*}}^{*}\right)=0$.
Then, we may identify such an optimal hypothesis index with the optimisation
problem
\[
k^{*}=\min\underset{k\in\left[m\right]}{\arg\min}\left\{ \min_{\bm{\theta}_{k}\in\mathbb{T}_{k}}r_{k}\left(\bm{\theta}_{k}\right)\right\} \text{,}
\]
which we can then estimate via the IC
\begin{equation}
\hat{K}_{n}=\min\underset{k\in\left[m\right]}{\arg\min}\text{ }\left\{ \min_{\bm{\theta}_{k}\in\mathbb{T}_{k}}\text{ }\frac{1}{n}\sum_{i=1}^{n}\left\Vert \bm{X}_{i}-\bm{\theta}_{k}\bm{\theta}_{k}^{\top}\bm{X}_{i}\right\Vert _{2}^{2}+P_{k,n}\right\} \text{,}\label{eq: PCA IC}
\end{equation}
where we can take $P_{k,n}$ the same as (\ref{eq: SWIC/BIC}) to
satisfy B2 or B2$^{*}$.

Since (\ref{eq: trace version}) is quadratic, we can then compute
the derivative with respect to $\bm{\theta}_{k}$ to be $\partial\ell_{k}\left(\bm{x};\bm{\theta}_{k}\right)/\partial\bm{\theta}_{k}=2\bm{x}\bm{x}^{\top}\bm{\theta}_{k}$.
We can thus verify A2 by assuming that the fourth moment $\text{E}_{X}\left\Vert \bm{X}\right\Vert _{4}^{4}$
is finite, which then also verifies A1.

Noting the form of (\ref{eq: r_k for PCA}), we then conclude, by
\citet[Thm.  23.2]{Shalev-Shwartz2014}, that $r_{k}$ is minimised
at any
\begin{equation}
\bm{\theta}_{k}^{*}=\left[\begin{array}{ccc}
\bm{\upsilon}_{\Pi\left(1\right)} & \cdots & \bm{\upsilon}_{\Pi\left(k\right)}\end{array}\right]\text{,}\label{eq: solutions of PCA}
\end{equation}
where $\bm{\upsilon}_{1},\dots,\bm{\upsilon}_{k}\in\mathbb{R}^{m}$
are orthonormal eigenvectors corresponding to the first $k$ largest
eigenvalues of $\text{E}\left(\bm{X}\bm{X}^{\top}\right)$ and $\Pi$,
again, is a permutation of $\left[k\right]$. Thus, C1 does not hold
unless we localise $\mathbb{T}_{k}$ on a neighborhood of one of the
isolated matrices. C2 is trivial, and C3 and C4 can then be verified
when $\text{E}\left(\bm{X}\bm{X}^{\top}\right)$ is positive definite,
using the results from the Appendix. We note that conditions under
which the BIC is consistent for determining the optimal subspace dimension
of data has been previously considered in the works of \citet{Zhao1986},
\citet{Fujikoshi2016}, and \citet{Bai2018}.

\section{\label{sec:Numerical-experiments}Numerical experiments}

The discussions thus far have been concerned with the theoretical
performance of PanIC and BIC-like estimators when sample sizes approach
infinity. However, practical behaviours of such estimators vary when
sample sizes are small. As such, we find it useful to compare criteria
in small sample settings to examine their relative performances.

For a strictly increasing sequence of constants $\left(c_{k}\right)_{k\in\left[m\right]}$,
we will say that the (general) BIC is an estimator of form (\ref{eq: general IC}),
where we choose $P_{n,k}=c_{k}n^{-1}\log n/2$, for each $k\in\left[m\right]$.
Since the SWIC defined in Remark \ref{rem: SWIC} require a choices
for the parameters $\alpha$ and $\beta$, we find it practical to
couple the choice of $\alpha$ with the BIC, by choosing constants
$\alpha=\alpha\left(\beta,\nu\right)>0$, which are defined as the
value of $\alpha$ that makes the order $\beta$ SWIC and the BIC
equal, when the sample size is $\nu\in\mathbb{N}$:
\[
\alpha\left(\beta,\nu\right)=\frac{\log\nu}{2\sqrt{\nu\log_{+}^{\left(\beta\right)}\left(\nu\right)}}\text{.}
\]
To facilitate the comparison of the orders $\beta\in\left\{ 1,2\right\} $
SWIC with the BIC, we will consider the respective choices of $\alpha\left(\beta,\nu\right)$,
with $\nu\in\left\{ 10^{3},10^{4}\right\} $. We will also provide
comparisons with the AIC ($P_{n,k}=c_{k}/n$), for completeness. 

We conduct all experiments in the $\mathsf{R}$ programming language
and make our code public at: \href{https://github.com/hiendn/PanIC}{https://github.com/hiendn/PanIC}.

\subsection{\label{subsec:Normal-finite-mixture}Normal finite mixture models}

In this experiment, we simulate data $\left(\bm{X}_{i}\right)_{i\in\left[n\right]}$
from bivariate normal mixture models with densities of forms
\[
\sum_{z=1}^{k^{*}}\pi_{z}\phi\left(\bm{x};\bm{\mu}_{z},\bm{\Sigma}_{z}\right)\text{,}
\]
where $\phi\left(\cdot;\bm{\mu},\bm{\Sigma}\right)$ is the normal
density with mean vector $\bm{\mu}\in\mathbb{R}^{2}$ and covariance
matrix $\bm{\Sigma}\in\mathbb{R}^{2\times2}$. Here, we consider four
different scenarios:
\begin{lyxlist}{00.00.0000}
\item [{S1.1}] $k^{*}=4$ and
\[
\left(\bm{\mu}_{z}\right)_{z\in\left[4\right]}=2\times\left\{ \left[\begin{array}{c}
1\\
1
\end{array}\right],\left[\begin{array}{c}
1\\
2
\end{array}\right],\left[\begin{array}{c}
2\\
1
\end{array}\right],\left[\begin{array}{c}
2\\
2
\end{array}\right]\right\} \text{.}
\]
\item [{S1.2}] $k^{*}=4$ and 
\[
\left(\bm{\mu}_{z}\right)_{z\in\left[4\right]}=3\times\left\{ \left[\begin{array}{c}
1\\
1
\end{array}\right],\left[\begin{array}{c}
1\\
2
\end{array}\right],\left[\begin{array}{c}
2\\
1
\end{array}\right],\left[\begin{array}{c}
2\\
2
\end{array}\right]\right\} \text{.}
\]
\item [{S1.3}] $k^{*}=6$ and 
\[
\left(\bm{\mu}_{z}\right)_{z\in\left[6\right]}=2\times\left\{ \left[\begin{array}{c}
1\\
1
\end{array}\right],\left[\begin{array}{c}
1\\
2
\end{array}\right],\left[\begin{array}{c}
1\\
3
\end{array}\right],\left[\begin{array}{c}
2\\
1
\end{array}\right],\left[\begin{array}{c}
2\\
2
\end{array}\right],\left[\begin{array}{c}
2\\
3
\end{array}\right]\right\} \text{.}
\]
\item [{S1.4}] $k^{*}=6$ and 
\[
\left(\bm{\mu}_{z}\right)_{z\in\left[6\right]}=3\times\left\{ \left[\begin{array}{c}
1\\
1
\end{array}\right],\left[\begin{array}{c}
1\\
2
\end{array}\right],\left[\begin{array}{c}
1\\
3
\end{array}\right],\left[\begin{array}{c}
2\\
1
\end{array}\right],\left[\begin{array}{c}
2\\
2
\end{array}\right],\left[\begin{array}{c}
2\\
3
\end{array}\right]\right\} \text{.}
\]
\end{lyxlist}
In all scenarios S1.1--S1.4, we set $\pi_{z}=1/k^{*}$ and $\bm{\Sigma}_{z}=\mathbf{I}_{2}$,
for each $z\in\left[k^{*}\right]$.

We take a sequence of hypothesis of form (\ref{eq: mixture hypotheses})
with $m=8$ and assess the performances of the AIC, BIC and SWIC of
order $\beta\in\left\{ 1,2\right\} $, using penalities of forms $k\left(q+1\right)/n$,
(\ref{eq: BIC mixture}), and (\ref{eq: SWIC mixture}), respectively,
with $q=5$. To evaluate performances, we consider sample sizes $n\in\left\{ 10^{2},10^{3},10^{4},10^{5}\right\} $
and apply the IC on 100 simulation runs, and compute the average value
of $\hat{K}_{n}$ obtained via each IC (\emph{Avg}), as well as the
proportion of occasions when $\hat{K}_{n}=k^{*}$ (\emph{Prop}). The
experimental results are presented in Tables \ref{tab: normal mixture results}.
Per Remark \ref{Rem: Parametric vs nonparametric}, we note that the
order selection procedure is nontrivial since we are operating in
the well-specified parametric setting where the data generating process
of $\left(X_{i}\right)_{i\in\left[n\right]}$ corresponds to one of
the hypotheses in $\left(\mathscr{H}_{k}\right)_{k\in\left[m\right]}$.

\begin{sidewaystable}
\caption{\label{tab: normal mixture results}Results from 100 simulation runs
for each scenario of the normal finite mixture models experiment.
Underlined values highlight either that the corresponding IC has either
the closest Avg to $k^{*}$ or the greatest Prop value for the indicated
scenario. }

\centering{}%
\begin{tabular}{|cc|cc|cc|cc|cc|cc|cc|}
\hline 
 & \multicolumn{1}{c}{} &  & \multicolumn{1}{c}{} &  & \multicolumn{1}{c}{} & \multicolumn{8}{c|}{SWIC $\left(\beta,\nu\right)$}\tabularnewline
Scenario & \multicolumn{1}{c}{$n$} & \multicolumn{2}{c}{AIC} & \multicolumn{2}{c}{BIC} & \multicolumn{2}{c}{$\left(1,10^{3}\right)$} & \multicolumn{2}{c}{$\left(1,10^{4}\right)$} & \multicolumn{2}{c}{$\left(2,10^{3}\right)$} & \multicolumn{2}{c|}{$\left(2,10^{4}\right)$}\tabularnewline
 & \multicolumn{1}{c}{} & Avg & \multicolumn{1}{c}{Prop} & Avg & \multicolumn{1}{c}{Prop} & Avg & \multicolumn{1}{c}{Prop} & Avg & \multicolumn{1}{c}{Prop} & Avg & \multicolumn{1}{c}{Prop} & Avg & Prop\tabularnewline
\hline 
\hline 
S1.1 & $10^{2}$ & 1.54 & 0.02 & 1.00 & 0.00 & 1.77 & 0.03 & 6.48 & \uline{0.08} & 1.55 & 0.02 & \uline{6.12} & \uline{0.08}\tabularnewline
$\left(k^{*}=4\right)$ & $10^{3}$ & \uline{2.28} & \uline{0.08} & 1.03 & 0.00 & 1.03 & 0.00 & 1.92 & 0.02 & 1.03 & 0.00 & 1.83 & 0.01\tabularnewline
 & $10^{4}$ & \uline{4.02} & \uline{0.88} & 3.32 & 0.32 & 2.19 & 0.00 & 3.32 & 0.32 & 2.36 & 0.00 & 3.32 & 3.32\tabularnewline
 & $10^{5}$ & 5.02 & 0.28 & 4.31 & 0.71 & 3.39 & 0.39 & \uline{4.07} & \uline{0.89} & 3.63 & 0.63 & 4.09 & 0.87\tabularnewline
\hline 
\hline 
S1.2 & $10^{2}$ & 2.66 & 0.09 & 1.06 & 0.00 & \uline{2.97} & \uline{0.14} & 6.86 & 0.10 & 2.72 & 0.10 & 6.52 & 0.11\tabularnewline
$\left(k^{*}=4\right)$ & $10^{3}$ & 4.05 & 0.97 & 3.53 & 0.54 & 3.53 & 0.54 & \uline{4.00} & \uline{1.00} & 3.53 & 0.54 & \uline{4.00} & \uline{1.00}\tabularnewline
 & $10^{4}$ & 4.00 & 1.00 & 4.00 & 1.00 & 4.00 & 1.00 & 4.00 & 1.00 & 4.00 & 1.00 & 4.00 & 1.00\tabularnewline
 & $10^{5}$ & 4.00 & 1.00 & 4.00 & 1.00 & 4.00 & 1.00 & 4.00 & 1.00 & 4.00 & 1.00 & 4.00 & 1.00\tabularnewline
\hline 
\hline 
S1.3 & $10^{2}$ & 1.45 & 0.01 & 1.00 & 0.00 & 1.69 & 0.04 & 6.84 & \uline{0.14} & 1.47 & 0.01 & \uline{6.00} & 0.10\tabularnewline
$\left(k^{*}=6\right)$ & $10^{3}$ & \uline{2.73} & 0.00 & 1.67 & 0.00 & 1.67 & 0.00 & 2.57 & 0.00 & 1.67 & 0.00 & 2.53 & 0.00\tabularnewline
 & $10^{4}$ & \uline{4.91} & \uline{0.20} & 3.84 & 0.00 & 2.37 & 0.00 & 3.84 & 0.00 & 2.44 & 0.00 & 3.84 & 0.00\tabularnewline
 & $10^{5}$ & 7.04 & 0.20 & \uline{5.89} & \uline{0.45} & 3.96 & 0.00 & 4.64 & 0.11 & 4.01 & 0.00 & 4.67 & 0.12\tabularnewline
\hline 
\hline 
S1.4 & $10^{2}$ & 2.84 & 0.04 & 1.07 & 0.00 & 3.30 & 0.06 & 6.59 & 0.21 & 3.02 & 0.05 & \uline{6.26} & \uline{0.23}\tabularnewline
$\left(k^{*}=6\right)$ & $10^{3}$ & \uline{5.64} & \uline{0.53} & 3.61 & 0.01 & 3.61 & 0.01 & 5.30 & 0.43 & 3.61 & 0.01 & 5.20 & 0.40\tabularnewline
 & $10^{4}$ & 6.11 & 0.90 & \uline{6.09} & \uline{0.91} & 5.75 & 0.84 & \uline{6.09} & \uline{0.91} & 5.87 & 0.89 & \uline{6.09} & \uline{0.91}\tabularnewline
 & $10^{5}$ & 6.18 & 0.85 & 6.14 & 0.87 & \uline{6.05} & \uline{0.95} & \uline{6.05} & \uline{0.95} & \uline{6.05} & \uline{0.95} & \uline{6.05} & \uline{0.95}\tabularnewline
\hline 
\end{tabular}
\end{sidewaystable}

We firstly notice that no IC is dominant across the range of scenarios
and sample sizes. However, we can make a number of generalisation.
Firstly, the choices of $\alpha$ constants for the SWIC leads to
criteria that are highly anti-conservative when $n=100$, in the sense
that the estimates $\hat{K}_{n}$ tend to be larger than the BIC and
even the AIC. This is typically a bad thing, since one generally would
like to choose the most parsimonious model, but we can observe that
this propensity to choose larger hypothesis classes can lead to greater
proportions of correct identifications of $k^{*}$.

This anti-conservativeness decreases dramatically, when comparing
the results for $n=1000$ to $n=100$, where we observe that the SWIC
become more conservative than the AIC, where the constants based on
$\nu=10^{4}$ tend to resemble the performance of the AIC in all scenarios. 

The AIC appears to perform well across all scenarios when $10^{4}$,
although, as expected by the fact that the AIC satisfies B1 but not
B2 (or B2$^{*}$), when $n$ is further increased to $10^{5}$, we
observe that the AIC becomes anti-conservative in all cases except
for S1.2, where all IC perform perfectly. Most importantly, we observe
that, as expected from Theorem \ref{Theorem: SWIC}, regardless of
the choice of $\alpha$ and $\beta$, the Avg values for the SWIC
appear to converge towards to $k^{*}$, as $n$ increases in all settings,
albeit rather slowly in the case of S1.3. We also observe the same
convergence regarding the BIC, which supports the theory of \citet{Leroux1992},
\citet{Keribin2000}, and \citet[Ch. 4]{Gassiat2018}.

\subsection{Least absolute deviation and $\epsilon$-support vector regression}

Let $\bm{X}=\left(Y,\bm{W}\right)\in\mathbb{R}\times\mathbb{R}^{m}$,
for $m=10$. We simulate $\left(\bm{X}_{i}\right)_{i\in\left[n\right]}$,
IID, from the same distribution as $\bm{X}$, defined via the linear
relationship
\[
Y=\bm{\theta}^{*\top}\bm{W}+E\text{,}
\]
where $\bm{W}\sim\text{Uniform}\left(\left[0,1\right]^{m}\right)$
and $E\sim\text{N}\left(0,1\right)$. Here, we consider four scenarios
S2.1--S2.4, whereupon in S2.1 and S2.2,
\[
\bm{\theta}^{*}=\left(1,0.75,0.5,0.25,0,0,0,0,0,0\right)^{\top}
\]
and in S2.3 and S2.4,
\[
\bm{\theta}^{*}=\left(1,0.85,0.7,0.55,0.4,0.25,0,0,0,0\right)^{\top}\text{.}
\]

Using the sequence of hypotheses of for(\ref{eq: SVR hypo}) with
$m=10$, we assess the performances of the AIC, BIC and SWIC of order
$\beta\in\left\{ 1,2\right\} $ of form (\ref{eq: SVR ICs}) with
penalties $k/n$, $kn^{-1}\log n/2$, and $\alpha k\sqrt{n^{-1}\log_{+}^{\left(\beta\right)}\left(n\right)}$,
respectively. Here, $\epsilon=0$ in S2.1 and S2.3, and $\epsilon=1$
in S2.2 and S2.4. We further note that $k^{*}=4$ in S2.1 and S2.2,
and $k^{*}=6$ in S2.3 and S2.4. To evaluate the performances of the
IC, we consider sample sizes $n\in\left\{ 10^{2},10^{3},10^{4},10^{5}\right\} $
and use 100 simulation runs of each scenario. Avg and Prop values
are recorded in Table \ref{tab: SVR}.

\begin{sidewaystable}
\caption{\label{tab: SVR}Results from 100 simulation runs for each scenario
of the LAD and SVR experiment. Underlined values highlight either
that the corresponding IC has either the closest Avg to $k^{*}$ or
the greatest Prop value for the indicated scenario.}

\centering{}%
\begin{tabular}{|cc|cc|cc|cc|cc|cc|cc|}
\hline 
 & \multicolumn{1}{c}{} &  & \multicolumn{1}{c}{} &  & \multicolumn{1}{c}{} & \multicolumn{8}{c|}{SWIC $\left(\beta,\nu\right)$}\tabularnewline
Scenario & \multicolumn{1}{c}{$n$} & \multicolumn{2}{c}{AIC} & \multicolumn{2}{c}{BIC} & \multicolumn{2}{c}{$\left(1,10^{3}\right)$} & \multicolumn{2}{c}{$\left(1,10^{4}\right)$} & \multicolumn{2}{c}{$\left(2,10^{3}\right)$} & \multicolumn{2}{c|}{$\left(2,10^{4}\right)$}\tabularnewline
 & \multicolumn{1}{c}{} & Avg & \multicolumn{1}{c}{Prop} & Avg & \multicolumn{1}{c}{Prop} & Avg & \multicolumn{1}{c}{Prop} & Avg & \multicolumn{1}{c}{Prop} & Avg & \multicolumn{1}{c}{Prop} & Avg & Prop\tabularnewline
\hline 
\hline 
S2.1 & $10^{2}$ & 3.62 & \uline{0.16} & 2.30 & 0.01 & \uline{4.18} & 0.14 & 7.97 & 0.07 & 3.66 & \uline{0.16} & 7.63 & 0.08\tabularnewline
$\left(k^{*}=4\right)$ & $10^{3}$ & 4.47 & 0.52 & 3.39 & 0.37 & 3.39 & 0.37 & 4.23 & 0.51 & 3.39 & 0.37 & \uline{4.12} & \uline{0.54}\tabularnewline
$\left(\epsilon=0\right)$ & $10^{4}$ & 5.13 & 0.60 & 4.02 & 0.98 & 3.98 & 0.98 & 4.02 & 0.98 & \uline{3.99} & \uline{0.99} & 4.02 & 0.98\tabularnewline
 & $10^{5}$ & 4.75 & 0.69 & \uline{4.00} & \uline{1.00} & \uline{4.00} & \uline{1.00} & \uline{4.00} & \uline{1.00} & \uline{4.00} & \uline{1.00} & \uline{4.00} & \uline{1.00}\tabularnewline
\hline 
\hline 
S2.2 & $10^{2}$ & 2.69 & 0.05 & 1.95 & 0.00 & \uline{2.76} & 0.06 & 6.06 & 0.13 & 2.70 & 0.05 & 5.34 & \uline{0.16}\tabularnewline
$\left(k^{*}=4\right)$ & $10^{3}$ & \uline{3.71} & \uline{0.61} & 3.13 & 0.17 & 3.13 & 0.17 & 3.59 & 0.57 & 3.13 & 0.17 & 3.55 & 0.55\tabularnewline
$\left(\epsilon=1\right)$ & $10^{4}$ & 4.17 & 0.89 & \uline{3.99} & \uline{0.99} & 3.79 & 0.79 & \uline{3.99} & \uline{0.99} & 3.87 & 0.87 & \uline{3.99} & \uline{0.99}\tabularnewline
 & $10^{5}$ & 4.22 & 0.92 & \uline{4.00} & \uline{1.00} & \uline{4.00} & \uline{1.00} & \uline{4.00} & \uline{1.00} & \uline{4.00} & \uline{1.00} & \uline{4.00} & \uline{1.00}\tabularnewline
\hline 
\hline 
S2.3 & $10^{2}$ & 4.94 & 0.14 & 3.45 & 0.01 & \uline{5.24} & \uline{0.17} & 7.79 & 0.09 & 5.08 & 0.16 & 7.56 & 0.08\tabularnewline
$\left(k^{*}=6\right)$ & $10^{3}$ & 6.76 & 0.55 & 5.37 & 0.38 & 5.37 & 0.38 & 6.39 & 0.58 & 5.37 & 0.38 & \uline{6.31} & \uline{0.60}\tabularnewline
$\left(\epsilon=0\right)$ & $10^{4}$ & 6.74 & 0.64 & \uline{6.02} & \uline{0.98} & 5.96 & 0.96 & \uline{6.02} & \uline{0.98} & 5.97 & 0.97 & \uline{6.02} & \uline{0.98}\tabularnewline
 & $10^{5}$ & 6.79 & 0.65 & 6.02 & 0.98 & \uline{6.00} & \uline{1.00} & \uline{6.00} & \uline{1.00} & \uline{6.00} & \uline{1.00} & \uline{6.00} & \uline{1.00}\tabularnewline
\hline 
\hline 
S2.4 & $10^{2}$ & 3.91 & 0.09 & 2.92 & 0.00 & 4.20 & 0.08 & 6.78 & 0.13 & 3.95 & 0.09 & \uline{6.29} & \uline{0.15}\tabularnewline
$\left(k^{*}=6\right)$ & $10^{3}$ & \uline{5.64} & \uline{0.50} & 4.89 & 0.12 & 4.89 & 0.12 & 5.42 & 0.41 & 4.89 & 0.12 & 5.36 & 0.38\tabularnewline
$\left(\epsilon=1\right)$ & $10^{4}$ & 6.14 & 0.91 & \uline{6.00} & \uline{1.00} & 5.64 & 0.64 & \uline{6.00} & \uline{1.00} & 5.80 & 0.80 & \uline{6.00} & \uline{1.00}\tabularnewline
 & $10^{5}$ & 6.10 & 0.91 & \uline{6.00} & \uline{1.00} & \uline{6.00} & \uline{1.00} & \uline{6.00} & \uline{1.00} & \uline{6.00} & \uline{1.00} & \uline{6.00} & \uline{1.00}\tabularnewline
\hline 
\end{tabular}
\end{sidewaystable}

As in Section \ref{subsec:Normal-finite-mixture}, we observe that
when $n=100$, the SWIC with the considered choices for $\alpha$
can be highly anti-conservative, especially when $\nu=10000$. This
can lead to over-estimation of $k^{*}$, although this can sometimes
lead to a higher value of Prop. For $n>100$, we observe that the
SWIC and BIC all tend to perform similarly, and in all cases the Avg
and Prop appear to converge, to $k^{*}$ and $1$, respectively, as
$n$ increases. In the case of the SWIC, this is predictable from
Theorem \ref{Theorem: SWIC}, and the consistency of the BIC is in
the LAD scenarios (S2.1 and S2.3) is also supported by the works of
\citet{Bai1998Estimation-of-m} and \citet{Lee2014Model-selection}.
However, the apparent consistency of the BIC in $\epsilon>0$ SVR
scenarios is of interest and can likely be established via the same
analysis as that of \citet{Lee2014Model-selection}. Proving this
result is outside of the scope the current work. We finally note that
the behaviours of the SWIC and BIC contrast with that of the AIC,
which appears to remain anti-conservative, even when $n$ is large,
in all scenarios.

\subsection{Principal component analysis}

We now simulate $\left(\bm{X}_{i}\right)_{i\in\left[n\right]}$, IID,
with the same distribution as $\bm{X}\in\mathbb{R}^{m}$ with $m=10$,
where
\[
\bm{X}=\mathbf{A}\bm{Y}
\]
for some fixed matrix $\mathbf{A}\in\mathbb{R}^{m\times k^{*}}$ and
$k^{*}$-dimensional random variable $\bm{Y}$ that differs in distribution
according to the simulation scenario. Here, we take $\bm{Y}\sim\text{N}\left(0,\mathbf{R}\right)$
in Scenarios S3.1 and S3.3, and $\bm{Y}\sim t_{5}\left(0,\mathbf{R}\right)$
in S3.2 and S3.4, where $\mathbf{R}$ has diagonal elements equal
to 1 and off-diagonal elements equal to $3/4$, and where $t_{\text{df}}$
denotes the multivariate Student-$t$ law with $\text{df}\in\mathbb{R}_{>0}$
degrees of freedom (cf. \citealt[Sec. 3.3.6]{Fang1990}). We consider
two different choices for $\mathbf{A}$:
\[
\left[\begin{array}{cccc}
1 & 0 & 0 & 0\\
0 & 0.25 & 0 & 0\\
0 & 0 & 0.1 & 0\\
0 & 0 & 0 & 0.1\\
1 & 0.25 & 0 & 0\\
0 & 0.25 & 0.1 & 0\\
0 & 0 & 0.1 & 0.1\\
1 & 0.25 & 0.1 & 0\\
0 & 0.25 & 0.1 & 0.1\\
1 & 0.25 & 0.1 & 0.1
\end{array}\right]\left(\text{S3.1, S3.2: }k^{*}=4\right)\text{, and}
\]

\[
\left[\begin{array}{cccccc}
1 & 0 & 0 & 0 & 0 & 0\\
0 & 0.5 & 0 & 0 & 0 & 0\\
0 & 0 & 0.25 & 0 & 0 & 0\\
0 & 0 & 0 & 0.25 & 0 & 0\\
0 & 0 & 0 & 0 & 0.1 & 0\\
0 & 0 & 0 & 0 & 0 & 0.1\\
1 & 0.5 & 0 & 0 & 0 & 0\\
0 & 0.5 & 0.25 & 0 & 0 & 0\\
0 & 0 & 0.25 & 0.25 & 0 & 0\\
0 & 0 & 0 & 0.25 & 0.1 & 0
\end{array}\right]\text{\ensuremath{\left(\text{S3.3, S3.4: }k^{*}=6\right)}.}
\]

Using the sequence of hypotheses of form (\ref{eq: SVR hypo}) with
$m=10$, we assess the performances of the AIC, BIC and SWIC of order
$\beta\in\left\{ 1,2\right\} $ of form (\ref{eq: SVR ICs}) with
penalties $km/n$, $kmn^{-1}\log n/2$, and $\alpha km\sqrt{n^{-1}\log_{+}^{\left(\beta\right)}\left(n\right)}$,
respectively. We again assess the IC via 100 simulation runs of each
scenario with sample sizes $n\in\left\{ 10^{2},10^{3},10^{4},10^{5}\right\} $,
and provide Avg and Prop values in Table (\ref{tab: PCA}).

\begin{sidewaystable}
\caption{\label{tab: PCA}Results from 100 simulation runs for each scenario
of the PCA experiment. Underlined values highlight either that the
corresponding IC has either the closest Avg to $k^{*}$ or the greatest
Prop value for the indicated scenario.}

\centering{}%
\begin{tabular}{|cc|cc|cc|cc|cc|cc|cc|}
\hline 
 & \multicolumn{1}{c}{} &  & \multicolumn{1}{c}{} &  & \multicolumn{1}{c}{} & \multicolumn{8}{c|}{SWIC $\left(\beta,\nu\right)$}\tabularnewline
Scenario & \multicolumn{1}{c}{$n$} & \multicolumn{2}{c}{AIC} & \multicolumn{2}{c}{BIC} & \multicolumn{2}{c}{$\left(1,10^{3}\right)$} & \multicolumn{2}{c}{$\left(1,10^{4}\right)$} & \multicolumn{2}{c}{$\left(2,10^{3}\right)$} & \multicolumn{2}{c|}{$\left(2,10^{4}\right)$}\tabularnewline
 & \multicolumn{1}{c}{} & Avg & \multicolumn{1}{c}{Prop} & Avg & \multicolumn{1}{c}{Prop} & Avg & \multicolumn{1}{c}{Prop} & Avg & \multicolumn{1}{c}{Prop} & Avg & \multicolumn{1}{c}{Prop} & Avg & Prop\tabularnewline
\hline 
\hline 
S3.1 & $10^{2}$ & 2.00 & 0.00 & 1.00 & 0.00 & 2.00 & 0.00 & \uline{2.99} & 0.00 & 2.00 & 0.00 & 2.90 & 0.00\tabularnewline
$\left(k^{*}=4\right)$ & $10^{3}$ & 3.00 & 0.00 & 3.00 & 0.00 & 3.00 & 0.00 & 3.00 & 0.00 & 3.00 & 0.00 & 3.00 & 0.00\tabularnewline
$\left(Y\sim\text{N}\right)$ & $10^{4}$ & \uline{4.00} & \uline{1.00} & \uline{4.00} & \uline{1.00} & 3.00 & 0.00 & \uline{4.00} & \uline{1.00} & 3.00 & 0.00 & \uline{4.00} & \uline{1.00}\tabularnewline
 & $10^{5}$ & 4.00 & 1.00 & 4.00 & 1.00 & 4.00 & 1.00 & 4.00 & 1.00 & 4.00 & 1.00 & 4.00 & 1.00\tabularnewline
\hline 
\hline 
S3.2 & $10^{2}$ & 2.11 & 0.00 & 1.72 & 0.00 & 2.24 & 0.00 & 3.00 & 0.00 & 2.14 & 0.00 & \uline{3.00} & 0.00\tabularnewline
$\left(k^{*}=4\right)$ & $10^{3}$ & \uline{3.24} & \uline{0.24} & 3.00 & 0.00 & 3.00 & 0.00 & 3.00 & 0.00 & 3.00 & 0.00 & 3.00 & 0.00\tabularnewline
$\left(Y\sim t_{5}\right)$ & $10^{4}$ & \uline{4.00} & \uline{1.00} & \uline{4.00} & \uline{1.00} & 3.00 & 0.00 & \uline{4.00} & \uline{1.00} & 3.00 & 0.00 & \uline{4.00} & \uline{1.00}\tabularnewline
 & $10^{5}$ & 4.00 & 1.00 & 4.00 & 1.00 & 4.00 & 1.00 & 4.00 & 1.00 & 4.00 & 1.00 & 4.00 & 1.00\tabularnewline
\hline 
\hline 
S3.3 & $10^{2}$ & 2.08 & 0.00 & 1.26 & 0.00 & 2.20 & 0.00 & \uline{3.10} & 0.00 & 2.12 & 0.00 & 3.00 & 0.00\tabularnewline
$\left(k^{*}=6\right)$ & $10^{3}$ & \uline{4.00} & 0.00 & 3.00 & 0.00 & 3.00 & 0.00 & \uline{4.00} & 0.00 & 3.00 & 0.00 & \uline{4.00} & 0.00\tabularnewline
$\left(Y\sim\text{N}\right)$ & $10^{4}$ & \uline{6.00} & \uline{1.00} & 4.97 & 0.00 & 4.00 & 0.00 & 4.97 & 0.00 & 4.00 & 0.00 & 4.97 & 0.00\tabularnewline
 & $10^{5}$ & \uline{6.00} & \uline{1.00} & \uline{6.00} & \uline{1.00} & 5.00 & 0.00 & \uline{6.00} & \uline{1.00} & 5.00 & 0.00 & \uline{6.00} & \uline{1.00}\tabularnewline
\hline 
\hline 
S3.4 & $10^{2}$ & 2.89 & 0.00 & 1.98 & 0.00 & 2.99 & 0.00 & \uline{3.89} & 0.00 & 2.93 & 0.00 & 3.65 & 0.00\tabularnewline
$\left(k^{*}=6\right)$ & $10^{3}$ & 4.00 & 0.00 & 4.00 & 0.00 & 4.00 & 0.00 & 4.00 & 0.00 & 4.00 & 0.00 & 4.00 & 0.00\tabularnewline
$\left(Y\sim t_{5}\right)$ & $10^{4}$ & \uline{6.00} & \uline{1.00} & 5.83 & 0.83 & 4.00 & 0.00 & 5.83 & 0.83 & 4.00 & 0.00 & 5.83 & 0.83\tabularnewline
 & $10^{5}$ & 6.00 & 1.00 & 6.00 & 1.00 & 6.00 & 1.00 & 6.00 & 1.00 & 6.00 & 1.00 & 6.00 & 1.00\tabularnewline
\hline 
\end{tabular}
\end{sidewaystable}

We again observe that the SWIC tend to be anti-conservative when $n=100$,
although we do not observe over-estimation of $k^{*}$ in any of our
scenarios, for any of the criteria. This is especially interesting
when inspecting the AIC results, are typically anti-conservative and
not consistent. When $n>100$ the SWIC and BIC appear to behave similarly,
with a tendency for Avg and Prop to approach $k^{*}$ and 1, respectively,
in accordance to Theorems \ref{Theorem: SWIC} and \ref{Theorem: BIC-like}.
We particularly note that $\nu=1000$ appears to be too conservative
when choosing $\alpha$, as in S2.3, both of the SWIC with $\nu=1000$
do not correctly identify $k^{*}$, even when $n=10^{5}$. This conservativeness
is further demonstrated when observing the results for $n=10^{4}$
in S3.1, S3.2 and S3.4, where the BIC appears to estimate $k^{*}$
correctly with high probability, but where the SWIC with $\nu=1000$
yield under-estimates of $k^{*}$. 

\section{Discussion}

Model selection is a ubiquitous task that arises when using many methods
for statistical inference and machine learning. The method of IC has
long been an essential tool for model selection in the likelihood
setting and related paradigms, such as quasi-likelihood and composite
likelihood-based methods. Using asymptotic theory from the stochastic
programming literature, we have demonstrated that simple and practical
IC-based methods are available for conducting variable selection in
general settings, such as when there are no clear definitions of likelihood
functions, and when estimators may not be unique. Following the primary
work of \citet{Sin:1996aa}, we propose the PanIC framework for model
selection, that requires minimal assumptions for guaranteeing consistent
identification of parsimonious models. We also provide relaxations
of sufficient conditions for proving consistency of BIC-like criteria,
when compared to \citet{Sin:1996aa} and \citet{Baudry:2015aa}. Our
empirical studies demonstrate that PanIC estimators can be usefully
applied in typical model selection settings where the consistency
of the BIC-like criteria are hard to verify, or when there are no
IC-based methods currently available.

Although general and useful, we have identified some directions to
improve upon the present work. Firstly, while restrictive, sufficient
conditions from \citet{Sin:1996aa} permit the possibility of non-IID
data, when verifying the consistency of PanIC and BIC-like criteria,
via generic non-IID strong laws of large numbers. Although non-trivial,
inspections of the proofs of Lemmas \ref{Lem: A1-A3} and \ref{Lem: A1-A3 C1--C5}
from \citet{Shapiro2021} reveal a reliance on empirical processes
methods for proving Donsker-type central limit theorems that can be
replaced by non-IID variants, such as those of \citet{Dedecker:2002aa}
and \citet[Ch. 8]{Rio:2017aa}. Progress in this direction has been
made by \citet{Wang:2022aa}, although their results are insufficient
for establishing non-IID versions of Lemmas \ref{Lem: A1-A3} and
\ref{Lem: A1-A3 C1--C5}.

Secondly, one may contend that our choices of $\alpha$ for the SWIC
in Section \ref{sec:Numerical-experiments} are arbitrary and unlikely
to be optimal, although we note that in the context of Theorem (\ref{Theorem: SWIC}),
the choice of $\alpha$ is immaterial when considering the asymptotic
properties of the SWIC and other PanIC. However, in more critical
applied work, it is possible to select $\alpha$ in a principled and
discerning manner. One such approach would be to conduct pilot simulations
of potential models from which data may arise and then apply the so-called
\emph{slope heuristics} of \citet{BirgeMassart2007} and \citet{Arlot:2020aa}
to estimate the value of $\alpha$. Such an approach has been usefully
applied, for example, in \citet{Nguyen:2018ac} and \citet{Nguyen:2022aa},
via the software of \citet{BaudryMaugisMichel2012}.

Thirdly, we note that all of our results are asymptotic and thus provide
no guarantees for any fixed sample size $n\in\mathbb{N}$. Finite
sample oracle inequalities for estimators of form (\ref{eq: general IC})
have been comprehensively studied in the works of \citet{Massart2007}
and \citet{Giraud2022}, for example. Although desirable, the establishment
of oracle inequalities may be excessive and require more bespoke and
precise analyses that go beyond the simple verification of the sufficient
conditions required by PanIC. We note that when $\left(\mathscr{H}_{k}\right)_{k\in\left[m\right]}$
is a sequence of nested hypotheses and when $\ell$ corresponds to
a likelihood, composite likelihood, or conditional likelihood maximisation
problem, then the approach of \citet{Nguyen:2022ab} provides a straightforward
method for obtaining finite sample confidence statements regarding
the event $\left\{ \omega:\tilde{K}_{n}\ge k^{*}\right\} $, where
$\tilde{K}_{n}$ is an estimator of $k^{*}$ that is obtained via
the \emph{universal inference} hypothesis testing framework of \citet{Wasserman:2020aa}.
Such confidence statements can be used in combination with PanIC to
provide complementary measures of uncertainty regarding the value
of $k^{*}$.

Lastly, although we have only concerned ourselves with sufficient
conditions that guarantee (\ref{eq: Consistency}), some situations
may require the guarantee of \emph{strong consistency}: 
\begin{equation}
\text{Pr}\left(\lim_{n\rightarrow\infty}\hat{K}_{n}=k^{*}\right)=1\text{.}\label{eq: strong consistency}
\end{equation}
Demanding regularity conditions are provided in \citet{Sin:1996aa}
for proving (\ref{eq: strong consistency}) in the general setting,
while setting-specific requirements are provided in the works of \citet{Zhao1986},
\citet{Keribin2000}, and \citet[Ch. 4]{Gassiat2018}, for example.
This prompts the question as to whether there is a version of Theorem
\ref{Theorem: SWIC}, with comparable assumptions, that guarantees
(\ref{eq: strong consistency}). We are happy to provide an affirmative
answer, thanks to the recent work of \citet{Banholzer:2022aa}. The
following result can be inferred from \citet[Thm. 4]{Banholzer:2022aa}.
\begin{lem}
\label{Lem: A1-A3 a.s.}Assume that $\left(X_{i}\right)_{i\in\left[n\right]}$
is an IID sequence, and that A1 and A2 hold, for each $k\in\left[m\right]$.
Then, for each $k$,
\begin{equation}
\sqrt{\frac{n}{\log_{+}^{\left(2\right)}\left(n\right)}}\left(\min_{\bm{\theta}_{k}\in\mathbb{T}_{k}}R_{k,n}\left(\bm{\theta}_{k}\right)-\min_{\bm{\theta}_{k}\in\mathbb{T}_{k}}r_{k}\left(\bm{\theta}_{k}\right)\right)=O_{\text{a.s.}}\left(1\right)\text{.}\label{eq: Lemma 2 a.s.}
\end{equation}
\end{lem}
Make the following assumptions, for each $k\in\left[m\right]$:
\begin{description}
\item [{B1$^{*}$}] $P_{k,n}>0$, for each $k\in\left[m\right]$ and $n\in\mathbb{N}$,
and $P_{k,n}=o_{\text{a.s.}}\left(1\right)$, as $n\rightarrow\infty$.
\item [{B2$^{**}$}] If $k<l$, then $\sqrt{n/\log_{+}^{\left(2\right)}\left(n\right)}\left\{ P_{l,n}-P_{k,n}\right\} \rightarrow\infty$,
almost surely, as $n\rightarrow\infty$.
\end{description}
\begin{thm}
\label{Theorem: SWIC a.s.}Assume that $\left(\bm{X}_{i}\right)_{i\in\left[n\right]}$
is an IID sequence, and that A1, A2, B1$^{*}$ and B2$^{**}$ hold,
for each $k\in\left[m\right]$. Then, the estimator (\ref{eq: general IC})
satisfies the strong consistency property (\ref{eq: strong consistency}).
\end{thm}
Theorem \ref{Theorem: SWIC a.s.} can be considered as a strong consistency
version of Theorem \ref{Theorem: SWIC}. We may refer to the class
of IC that satisfies B1$^{*}$ and B2$^{**}$ as \emph{strong PanIC}.
Note that SWIC of order 1 are strong PanIC but SWIC with $\beta>1$
and the BIC are not. This provides some incentive towards the use
of the order 1 PanIC over those of higher order, although we observe
from Section \ref{sec:Numerical-experiments} that the difference
may not be so pronounced in finite samples. 

We note that the BIC and Hannan--Quinn IC can be proved to be strongly
consistent in many settings, but do not satisfy the assumptions of
Theorem \ref{Theorem: SWIC a.s.}. Unfortunately, we cannot yet prove
a strong consistency equivalent to Theorem \ref{Theorem: BIC-like}
that would apply to the BIC and Hannan--Quinn IC. We believe that
such a result is available and leave its proof to future research.

\section*{Appendix}

\subsection*{Technical definitions and results}

\subsubsection*{Assumptions C3 and C4}

The following definitions follow the exposition of \citet[Sec. 9.1.5]{Shapiro2021}.
Let $\mathbb{T}\subset\mathbb{R}^{q}$ for some $q\in\mathbb{N}$.
We firstly define
\[
\mathcal{T}_{\mathbb{T}}^{2}\left(\bm{\theta},\bm{\eta}\right)=\left\{ \bm{\tau}\in\mathbb{R}^{q}:\text{dist}\left(\bm{\theta}+t\bm{\eta}+\frac{1}{2}t^{2}\bm{\tau},\mathbb{T}\right)=o\left(t^{2}\right),t\ge0\right\} 
\]
to be the \emph{second order tangent set} of $\mathbb{T}$, at the
point $\bm{\theta}\in\mathbb{T}$, in the direction $\bm{\eta}\in\mathcal{T}_{\mathbb{T}}\left(\bm{\theta}\right)$.
Here, $\mathcal{T}_{\mathbb{T}}\left(\bm{\theta}\right)$ is the so-called
the \emph{contingent cone} to $\mathbb{T}$ at $\bm{\theta}\in\mathbb{T}$,
defined as the set of vectors $\bm{\eta}\in\mathbb{R}^{q}$, such
that there exist a sequences $\left(\bm{\eta}_{j}\right)_{j\in\mathbb{N}}\subset\mathbb{R}^{q}$
and $\left(t_{j}\right)_{j\in\mathbb{N}}\subset\mathbb{R}_{\ge0}$,
where $\lim_{j\rightarrow\infty}\bm{\eta}_{j}=\bm{\eta}$ and $\bm{\theta}+t_{j}\bm{\eta}_{j}\in\mathbb{T}$,
for each $j\in\mathbb{N}$. We then say that $\mathbb{T}$ is \emph{second
order regular} at $\bm{\theta}^{*}\in\mathbb{T}$ if for any sequence
$\left(\bm{\theta}_{j}\right)_{j\in\mathbb{N}}\subset\mathbb{T}$
of the form 
\[
\bm{\theta}_{j}=\bm{\theta}^{*}+t_{j}\bm{\eta}+\frac{1}{2}t_{j}^{2}\bm{\tau}_{j}\text{,}
\]
where $\lim_{j\rightarrow\infty}t_{j}=0$ and $\lim_{j\rightarrow\infty}t_{j}\bm{\tau}_{j}=0$
for sequences $\left(t_{j}\right)_{j\in\mathbb{N}}\subset\mathbb{R}_{\ge0}$
and $\left(\bm{\tau}_{j}\right)_{j\in\mathbb{N}}\subset\mathbb{R}^{q}$,
it follows that
\[
\lim_{j\rightarrow\infty}\text{dist}\left(\bm{\tau}_{j},\mathcal{T}_{\mathbb{T}}^{2}\left(\bm{\theta}^{*},\bm{\eta}\right)\right)=0\text{.}
\]

Next, let $r:\mathbb{T}\rightarrow\mathbb{R}$ be a function of interest.
We say that $r$ satisfies the \emph{quadratic growth condition} at
$\bm{\theta}^{*}\in\mathbb{T}$ if there exists a constant $\mathfrak{C}>0$
and a neighborhood $\mathcal{N}\subset\mathbb{R}^{q}$ of $\bm{\theta}^{*}$,
such that
\[
r\left(\bm{\theta}\right)\ge r\left(\bm{\theta}^{*}\right)+\mathfrak{C}\left\Vert \bm{\theta}-\bm{\theta}^{*}\right\Vert ^{2}\text{,}
\]
for all $\bm{\theta}\in\mathcal{N}\cap\mathbb{T}$. 

As a consequence of \citet[Prop. 3.88]{Bonnans2000Perturbation-An},
we have the fact that if 
\[
\mathbb{T}=\left\{ \bm{\theta}\in\mathbb{R}^{q}:g_{u}\left(\bm{\theta}\right)\le0\text{, }u\in\left[\mathfrak{m}\right],h_{v}\left(\bm{\theta}\right)=0\text{, }v\in\left[\mathfrak{n}\right]\right\} 
\]
for some $\mathfrak{m},\mathfrak{n}\in\mathbb{N}$, where each $g_{u},h_{v}:\mathbb{R}^{q}\rightarrow\mathbb{R}$
is twice continuously differentiable, and the point $\bm{\theta}^{*}$
satisfies the so-called \emph{Mangasarian--Fromovitz constraint qualification},
then $\mathbb{T}$ is second order regular at $\bm{\theta}^{*}$.
Here, the Mangasarian--Fromovitz constraint qualification can be
stated as follows: the vectors $\left.\partial h_{v}\left(\bm{\theta}\right)/\partial\bm{\theta}\right|_{\bm{\theta}=\bm{\theta}^{*}}$
($v\in\left[\mathfrak{n}\right]$) are linearly independent, and there
exists a vector $\bm{\lambda}\in\mathbb{R}^{q}$, such that $\lambda^{\top}\left.\partial h_{v}\left(\bm{\theta}\right)/\partial\bm{\theta}\right|_{\bm{\theta}=\bm{\theta}^{*}}=0$,
for each $v\in\left[n\right]$, and $\bm{\lambda}^{\top}\left.\partial g_{u}\left(\bm{\theta}\right)/\partial\bm{\theta}\right|_{\bm{\theta}=\bm{\theta}^{*}}<0$,
for each $u\in\left[\mathfrak{m}\right]$, such that $g_{u}\left(\bm{\theta}^{*}\right)=0$.
In particular, $\bm{\theta}^{*}\in\mathbb{T}$ is always second order
regular if $\mathbb{T}$ is \emph{polyhedral}, in the sense that $g_{u}$
and $h_{v}$ are affine functions of $\bm{\theta}\in\mathbb{T}$,
for each $u\in\left[\mathfrak{m}\right]$ and $v\in\left[\mathfrak{n}\right]$.

Using \citet[Props. 9.16 and 9.21]{Shapiro2021}, we have the facts
that if $\bm{\theta}^{*}\in\mathbb{T}$ is a global minimum of $r$,
and $r$ is continuously differentiable at $\bm{\theta}^{*}$, then
$\bm{\eta}^{\top}\left.\partial r\left(\bm{\theta}\right)/\partial\bm{\theta}\right|_{\bm{\theta}=\bm{\theta}^{*}}\ge0$,
for every $\bm{\eta}\in\mathcal{T}_{\mathbb{T}}\left(\bm{\theta}^{*}\right)$,
and if $\mathbb{T}$ is second order regular and $r$ is twice continuously
differentiable at $\bm{\theta}^{*}$, then $r$ satisfies the quadratic
growth condition at $\bm{\theta}^{*}$, if and only if,
\begin{equation}
\bm{\eta}^{\top}\left.\frac{\partial^{2}r\left(\bm{\theta}\right)}{\partial\bm{\theta}\partial\bm{\theta}^{\top}}\right|_{\bm{\theta}=\bm{\theta}^{*}}\bm{\eta}+\inf_{\bm{\lambda}\in\mathcal{T}_{\mathbb{T}}^{2}\left(\bm{\theta}^{*},\bm{\eta}\right)}\bm{\lambda}^{\top}\left.\frac{\partial r\left(\bm{\theta}\right)}{\partial\bm{\theta}}\right|_{\bm{\theta}=\bm{\theta}^{*}}>0\text{,}\label{eq: nec/suf condition}
\end{equation}
for every $\bm{\eta}\in\mathcal{C}\left(\bm{\theta}^{*}\right)\backslash\left\{ \bm{0}\right\} $,
where
\[
\mathcal{C}\left(\bm{\theta}^{*}\right)=\left\{ \bm{\eta}\in\mathcal{T}_{\mathbb{T}}\left(\bm{\theta}^{*}\right):\bm{\eta}^{\top}\left.\frac{\partial r\left(\bm{\theta}\right)}{\partial\bm{\theta}}\right|_{\bm{\theta}=\bm{\theta}^{*}}=0\right\} 
\]
is the so-called \emph{critical cone} at the minimiser $\bm{\theta}^{*}$.
As noted in \citet{Shapiro2000}, the second term on the LHS of (\ref{eq: nec/suf condition})
is equal to zero if $\left.\partial r\left(\bm{\theta}\right)/\partial\bm{\theta}\right|_{\bm{\theta}=\bm{\theta}^{*}}=\bm{0}$
or if $\mathbb{T}$ is polyhedral.

\subsection*{Proofs}

\subsubsection*{Proof of Theorem \ref{Theorem: SWIC}}

Our proof largely follows the approach as that of \citet[Thm. 8.1]{Baudry:2015aa}.
Firstly, the IID assumption together with A1 and A2 permits the application
of a uniform strong law of large numbers (e.g., \citeauthor[Thm. 9.60]{Shapiro2021}),
which then implies that, for each $k$,
\begin{equation}
\min_{\bm{\theta}_{k}\in\mathbb{T}_{k}}R_{k,n}\left(\bm{\theta}_{k}\right)\stackrel[n\rightarrow\infty]{\text{a.s.}}{\longrightarrow}\min_{\bm{\theta}_{k}\in\mathbb{T}_{k}}r_{k}\left(\bm{\theta}_{k}\right)\text{,}\label{eq: SLLN min}
\end{equation}
by \citet[Prop. 5.2]{Shapiro2021}. 

Let
\[
\mathcal{K}=\underset{k\in\left[m\right]}{\arg\min}\left\{ \min_{\bm{\theta}_{k}\in\mathbb{T}_{k}}r_{k}\left(\bm{\theta}_{k}\right)\right\} 
\]
and for each $k\notin\mathcal{K}$, let 
\begin{equation}
\epsilon=\frac{1}{2}\left\{ \min_{\bm{\theta}_{k}\in\mathbb{T}_{k}}r_{k}\left(\bm{\theta}_{k}\right)-\min_{\bm{\theta}_{k^{*}}\in\mathbb{T}_{k^{*}}}r_{k^{*}}\left(\bm{\theta}_{k^{*}}\right)\right\} >0\text{.}\label{eq: epsilon proof}
\end{equation}
Then, (\ref{eq: SLLN min}) and B1 imply that for every $\delta>0$,
\begin{equation}
\left|\min_{\bm{\theta}_{k}\in\mathbb{T}_{k}}R_{k,n}\left(\bm{\theta}_{k}\right)-\min_{\bm{\theta}_{k}\in\mathbb{T}_{k}}r_{k}\left(\bm{\theta}_{k}\right)\right|\le\frac{\epsilon}{3}\text{,}\label{eq: Difference in k}
\end{equation}
\begin{equation}
\left|\min_{\bm{\theta}_{k^{*}}\in\mathbb{T}_{k^{*}}}R_{k^{*},n}\left(\bm{\theta}_{k^{*}}\right)-\min_{\bm{\theta}_{k^{*}}\in\mathbb{T}_{k^{*}}}r_{k^{*}}\left(\bm{\theta}_{k^{*}}\right)\right|\le\frac{\epsilon}{3}\text{,}\label{eq: Difference in k*}
\end{equation}
and $P_{k,n}\le\epsilon/3$ all occur simultaneously with probability
at least $1-\delta$, whenever $n$ is sufficiently large (i.e., $n\ge n_{\delta}$,
for some large $n_{\delta}$). Thus,
\begin{align*}
\min_{\bm{\theta}_{k}\in\mathbb{T}_{k}}R_{k,n}\left(\bm{\theta}_{k}\right)+P_{k,n} & \ge\min_{\bm{\theta}_{k}\in\mathbb{T}_{k}}r_{k}\left(\bm{\theta}_{k}\right)-\frac{\epsilon}{3}\\
 & =\min_{\bm{\theta}_{k^{*}}\in\mathbb{T}_{k^{*}}}r_{k^{*}}\left(\bm{\theta}_{k^{*}}\right)+\frac{5\epsilon}{3}\\
 & \ge\min_{\bm{\theta}_{k^{*}}\in\mathbb{T}_{k^{*}}}R_{k^{*},n}\left(\bm{\theta}_{k^{*}}\right)+P_{k^{*},n}+\epsilon
\end{align*}
also occurs with probability at least $1-\delta$, for sufficiently
large $n$. But since $\epsilon>0$, we have the fact that 
\begin{equation}
\min_{\bm{\theta}_{k}\in\mathbb{T}_{k}}R_{k,n}\left(\bm{\theta}_{k}\right)+P_{k,n}>\min_{\bm{\theta}_{k^{*}}\in\mathbb{T}_{k^{*}}}R_{k^{*},n}\left(\bm{\theta}_{k^{*}}\right)+P_{k^{*},n}\text{,}\label{eq: smallest smaller than bigger}
\end{equation}
and thus by definition of (\ref{eq: general IC}), $\hat{K}_{n}\ne k$
occurs with probability at least $1-\delta$, for large $n$, for
each $k\notin\mathcal{K}$. Because there is only a finite number
of $k\notin\mathcal{K}$, we have the fact that 
\[
\lim_{n\rightarrow\infty}\text{Pr}\left(\hat{K}_{n}\notin\mathcal{K}\right)=0\text{.}
\]

Next, we consider when $k\in\mathcal{K}$ but $k>k^{*}$. From A1
and A2, Lemma \ref{eq: Lemma 1} implies that for each $k\in\mathcal{K}$,
\[
\sqrt{n}\left(\min_{\bm{\theta}_{k}\in\mathbb{T}_{k}}R_{k,n}\left(\bm{\theta}_{k}\right)-\min_{\bm{\theta}_{k}\in\mathbb{T}_{k}}r_{k}\left(\bm{\theta}_{k}\right)\right)=O_{\text{Pr}}\left(1\right)
\]
which implies that
\[
\sqrt{n}\left(\min_{\bm{\theta}_{k}\in\mathbb{T}_{k}}R_{k,n}\left(\bm{\theta}_{k}\right)-\min_{\bm{\theta}_{k^{*}}\in\mathbb{T}_{k^{*}}}R_{k^{*},n}\left(\bm{\theta}_{k^{*}}\right)\right)=O_{\text{Pr}}\left(1\right)\text{,}
\]
since $\min_{\bm{\theta}_{k^{*}}\in\mathbb{T}_{k^{*}}}r_{k^{*}}\left(\bm{\theta}_{k^{*}}\right)=\min_{\bm{\theta}_{k}\in\mathbb{T}_{k}}r_{k}\left(\bm{\theta}_{k}\right)$,
by definition of $\mathcal{K}$. Thus, for each $k$ and $\delta>0$,
there exists a constant $\mathfrak{M}>0$, such that
\[
\sqrt{n}\left|\min_{\bm{\theta}_{k}\in\mathbb{T}_{k}}R_{k,n}\left(\bm{\theta}_{k}\right)-\min_{\bm{\theta}_{k^{*}}\in\mathbb{T}_{k^{*}}}R_{k^{*},n}\left(\bm{\theta}_{k^{*}}\right)\right|\le\mathfrak{M}
\]
with probability at least $1-\delta$, for sufficiently large $n$.
On the other hand, by B2, for each $\delta>0$ and $\mathfrak{M}>0$,
\[
\sqrt{n}\left\{ P_{k,n}-P_{k^{*},n}\right\} >\mathfrak{M}
\]
with probability at least $1-\delta$, for sufficiently large $n$.
Thus, for sufficiently large $n$, with probability at least $1-2\delta$,
we have 
\[
\min_{\bm{\theta}_{k^{*}}\in\mathbb{T}_{k^{*}}}R_{k^{*},n}\left(\bm{\theta}_{k^{*}}\right)-\min_{\bm{\theta}_{k}\in\mathbb{T}_{k}}R_{k,n}\left(\bm{\theta}_{k}\right)\le\frac{\mathfrak{M}}{\sqrt{n}}<P_{k,n}-P_{k^{*},n}
\]
and thus
\[
\min_{\bm{\theta}_{k^{*}}\in\mathbb{T}_{k^{*}}}R_{k^{*},n}\left(\bm{\theta}_{k^{*}}\right)+P_{k^{*},n}<\min_{\bm{\theta}_{k}\in\mathbb{T}_{k}}R_{k,n}\left(\bm{\theta}_{k}\right)+P_{k,n}\text{,}
\]
for each $k\in\mathcal{K}$, such that $k>k^{*}$. Again, since there
is a finite number of elements of $\mathcal{K}$, and by definition
of (\ref{eq: general IC}), we have 
\[
\lim_{n\rightarrow\infty}\text{Pr}\left(\hat{K}_{n}\in\mathcal{K},\hat{K}_{n}>k^{*}\right)=0.
\]
Finally, we obtain the desired result via the limit:
\[
\lim_{n\rightarrow\infty}\text{Pr}\left(\hat{K}_{n}\ne k^{*}\right)=\lim_{n\rightarrow\infty}\left\{ \text{Pr}\left(\hat{K}_{n}\notin\mathcal{K}\right)+\text{Pr}\left(\hat{K}_{n}\in\mathcal{K},\hat{K}_{n}>k^{*}\right)\right\} =0\text{.}
\]

\subsubsection*{Proof of Theorem \ref{Theorem: BIC-like}}

Under A1, A2, and C1--C5, \citet[Thm. 4.4]{Shapiro2000} implies
the following result (via an application of \citealp[Thm. 9.56]{Shapiro2021}).
\begin{lem}
\label{Lem: A1-A3 C1--C5}Assume that $\left(\bm{X}_{i}\right)_{i\in\left[n\right]}$
is an IID sequence, and that A1, A2, and C1--C5 hold, for each $k\in\left[m\right]$.
Then, for each $k$,
\begin{equation}
n\left(\min_{\bm{\theta}_{k}\in\mathbb{T}_{k}}R_{k,n}\left(\bm{\theta}_{k}\right)-R_{k,n}\left(\bm{\theta}_{k}^{*}\right)\right)\rightsquigarrow\frac{1}{2}\varphi_{k}\left(Z_{k}\right)\text{,}\label{eq: Lemma 2}
\end{equation}
where $Z_{k}$ is normally distributed with mean zero and covariance
matrix
\[
\mathrm{E}\left\{ \left[\left.\frac{\partial\ell_{k}\left(X;\bm{\theta}_{k}\right)}{\partial\bm{\theta}_{k}}\right|_{\bm{\theta}_{k}=\bm{\theta}_{k}^{*}}\right]\left[\left.\frac{\partial\ell_{k}\left(X;\bm{\theta}_{k}\right)}{\partial\bm{\theta}_{k}}\right|_{\bm{\theta}_{k}=\bm{\theta}_{k}^{*}}\right]^{\top}\right\} -\left[\left.\frac{\partial r_{k}\left(\bm{\theta}\right)}{\partial\bm{\theta}_{k}}\right|_{\bm{\theta}_{k}=\bm{\theta}_{k}^{*}}\right]\left[\left.\frac{\partial r_{k}\left(\bm{\theta}\right)}{\partial\bm{\theta}_{k}}\right|_{\bm{\theta}_{k}=\bm{\theta}_{k}^{*}}\right]^{\top}
\]
and
\[
\varphi_{k}\left(\bm{\zeta}\right)=\inf_{\bm{\eta}\in\mathcal{C}\left(\bm{\theta}^{*}\right)}\left\{ 2\bm{\eta}^{\top}\bm{\zeta}+\bm{\eta}^{\top}\left.\frac{\partial^{2}r\left(\bm{\theta}\right)}{\partial\bm{\theta}_{k}\partial\bm{\theta}_{k}^{\top}}\right|_{\bm{\theta}_{k}=\bm{\theta}_{k}^{*}}\bm{\eta}+\inf_{\bm{\lambda}\in\mathcal{T}_{\mathbb{T}_{k}}^{2}\left(\bm{\theta}_{k}^{*},\bm{\eta}\right)}\bm{\lambda}^{\top}\left.\frac{\partial r_{k}\left(\bm{\theta}\right)}{\partial\bm{\theta}_{k}}\right|_{\bm{\theta}_{k}=\bm{\theta}_{k}^{*}}\right\} \text{.}
\]
\end{lem}
As per the use of Lemma \ref{Lem: A1-A3}, we only require Lemma \ref{Lem: A1-A3 C1--C5}
to obtain the fact that
\[
n\left(\min_{\bm{\theta}_{k}\in\mathbb{T}_{k}}R_{k,n}\left(\bm{\theta}_{k}\right)-R_{k,n}\left(\bm{\theta}_{k}^{*}\right)\right)=O_{\text{Pr}}\left(1\right)\text{,}
\]
for each $k\in\left[m\right]$. Recalling the definition of $\mathcal{K}$
from the previous proof, via C5 and Lemma \ref{Lem: A1-A3 C1--C5},
we then have
\begin{align*}
 & n\left(\min_{\bm{\theta}_{k}\in\mathbb{T}_{k}}R_{k,n}\left(\bm{\theta}_{k}\right)-\min_{\bm{\theta}_{k^{*}}\in\mathbb{T}_{k^{*}}}R_{k^{*},n}\left(\bm{\theta}_{k^{*}}\right)\right)\\
= & n\left(\min_{\bm{\theta}_{k}\in\mathbb{T}_{k}}R_{k,n}\left(\bm{\theta}_{k}\right)-R_{k,n}\left(\bm{\theta}_{k}^{*}\right)\right)-n\left(\min_{\bm{\theta}_{k^{*}}\in\mathbb{T}_{k^{*}}}R_{k^{*},n}\left(\bm{\theta}_{k^{*}}\right)-R_{k^{*},n}\left(\bm{\theta}_{k^{*}}^{*}\right)\right)\\
 & -n\left(R_{k,n}\left(\bm{\theta}_{k}^{*}\right)-R_{k^{*},n}\left(\bm{\theta}_{k^{*}}^{*}\right)\right)\\
= & O_{\text{Pr}}\left(1\right)+O_{\text{Pr}}\left(1\right)+O_{\text{Pr}}\left(1\right)=O_{\text{Pr}}\left(1\right)\text{,}
\end{align*}
for each $k\in\mathcal{K}$. The remainder of the proof follows in
the same way as that of Theorem \ref{Theorem: SWIC}. 

\subsubsection*{Proof of Theorem \ref{Theorem: SWIC a.s.}}

The proof follows in the same manner as that of Theorem \ref{Theorem: SWIC}.
By (\ref{eq: SLLN min}) and B1$^{*}$, for each $k\notin\mathcal{K}$,
we have both (\ref{eq: Difference in k}) and (\ref{eq: Difference in k*}),
almost surely when $n$ is sufficiently large, for $\epsilon$ defined
by (\ref{eq: epsilon proof}). This then implies (\ref{eq: smallest smaller than bigger}),
almost surely, for each $k\notin\mathcal{K}$, and thus, we can conclude
that 
\[
\text{Pr}\left(\lim_{n\rightarrow\infty}\hat{K}_{n}\notin\mathcal{K}\right)=0\text{,}
\]
since $\left[m\right]$ is finite.

Next, when A1 and A2 hold, and when $k\in\mathcal{K}$ and $k>k^{*}$,
Lemma \ref{Lem: A1-A3 a.s.} implies that 
\[
\sqrt{\frac{n}{\log_{+}^{\left(2\right)}\left(n\right)}}\left(\min_{\bm{\theta}_{k}\in\mathbb{T}_{k}}R_{k,n}\left(\bm{\theta}_{k}\right)-\min_{\bm{\theta}_{k}\in\mathbb{T}_{k}}r_{k}\left(\bm{\theta}_{k}\right)\right)=O_{\text{a.s.}}\left(1\right)
\]
and hence
\[
\sqrt{\frac{n}{\log_{+}^{\left(2\right)}\left(n\right)}}\left(\min_{\bm{\theta}_{k}\in\mathbb{T}_{k}}R_{k,n}\left(\bm{\theta}_{k}\right)-\min_{\bm{\theta}_{k^{*}}\in\mathbb{T}_{k^{*}}}R_{k^{*},n}\left(\bm{\theta}_{k^{*}}\right)\right)=O_{\text{a.s.}}\left(1\right)\text{,}
\]
by definition of $\mathcal{K}$. Thus, for each $k\in\mathcal{K}$,
there exists a constant $\mathfrak{M}>0$ such that
\begin{equation}
\sqrt{\frac{n}{\log_{+}^{\left(2\right)}\left(n\right)}}\left|\min_{\bm{\theta}_{k}\in\mathbb{T}_{k}}R_{k,n}\left(\bm{\theta}_{k}\right)-\min_{\bm{\theta}_{k^{*}}\in\mathbb{T}_{k^{*}}}R_{k^{*},n}\left(\bm{\theta}_{k^{*}}\right)\right|\le\mathfrak{M}\label{eq: bounded almost surely}
\end{equation}
for sufficiently large $n$, almost surely. Then, by B2$^{**}$, for
each $\mathfrak{M}>0$
\[
\sqrt{\frac{n}{\log_{+}^{\left(2\right)}\left(n\right)}}\left\{ P_{k,n}-P_{k^{*},n}\right\} >\mathfrak{M}
\]
for sufficiently large $n$, almost surely, which together with (\ref{eq: bounded almost surely}),
implies 
\[
\min_{\bm{\theta}_{k^{*}}\in\mathbb{T}_{k^{*}}}R_{k^{*},n}\left(\bm{\theta}_{k^{*}}\right)-\min_{\bm{\theta}_{k}\in\mathbb{T}_{k}}R_{k,n}\left(\bm{\theta}_{k}\right)\le\sqrt{\frac{\log_{+}^{\left(2\right)}\left(n\right)}{n}}\mathfrak{M}<P_{k,n}-P_{k^{*},n}
\]
and hence
\[
\min_{\bm{\theta}_{k^{*}}\in\mathbb{T}_{k^{*}}}R_{k^{*},n}\left(\bm{\theta}_{k^{*}}\right)+P_{k^{*},n}<\min_{\bm{\theta}_{k}\in\mathbb{T}_{k}}R_{k,n}\left(\bm{\theta}_{k}\right)+P_{k,n}\text{,}
\]
for sufficiently large $n$, almost surely, as required. Again, since
$\mathcal{K}$ is finite, the intersect of the required almost sure
events is also almost sure, and hence
\[
\text{Pr}\left(\lim_{n\rightarrow\infty}\hat{K_{n}}\in\mathcal{K},\lim_{n\rightarrow\infty}K_{n}>k^{*}\right)=0
\]
which then implies the desired result
\[
\text{Pr}\left(\lim_{n\rightarrow\infty}\hat{K}_{n}\ne k^{*}\right)=\text{Pr}\left(\lim_{n\rightarrow\infty}\hat{K}_{n}\notin\mathcal{K}\right)+\text{Pr}\left(\lim_{n\rightarrow\infty}\hat{K_{n}}\in\mathcal{K},\lim_{n\rightarrow\infty}K_{n}>k^{*}\right)=0\text{.}
\]

\bibliographystyle{plainnat}
\bibliography{2023_MathsETC}

\end{document}